\documentclass[11pt]{article}

\usepackage{amssymb,epsfig}
\usepackage{bbm}

\newcommand{\U}{\mathbb{U}}
\newcommand{\R}{\mathbb{R}}

\newcommand{\Z}{\mathbb{Z}}
\renewcommand{\P}{\mathbb{P}}
\newcommand{\Q}{\mathbb{Q}}
\newcommand{\E}{\mathbb{E}}

\newcommand{\N}{\mathbb{N}}

\newcommand{\T}{{\bf T}}

\newcommand{\B}{\mathbb{B}}
\renewcommand{\Z}{\mathbb{Z}}

\renewcommand{\S}{\mathbb{S}}

\def\build#1_#2^#3{\mathrel{
\mathop{\kern 0pt#1}\limits_{#2}^{#3}}}

\def\cq{$\hfill \square$}
\def \un{\underline}
\def\ind{{\mathbbm{1}}_}

\def\al{\alpha}
\def\vep{\varepsilon}

\def\vph{\varphi}
\def\lm{\lambda}
\def\sg{\sigma}
\def\dl{\delta}
\def\ka{\kappa}

\def\rad{{\mathcal R}}
\def\I{{\mathcal I}}
\def\dd{\partial}

\def\m{{\mathcal M}}

\def\g{{\mathcal G}}
\def\f{{\mathcal F}}

\def\u{{\mathcal U}}

\def\h{{\mathcal H}}
\def\t{{\mathcal T}}
\def\v{{\mathcal V}}

\def\s{{\mathcal S}}

\def\DD{{\mathcal D}}

\def\d{{\rm d}}

\def\lg{{\bf l}}
\def\e{{\bf e}}
\def\r{{\bf r}}
\def\q{{\bf q}}
\def\tr{{\bf t}}
\def\uu{{\bf u}}

\def\mm{{\bf m}}
\def\pp{{\bf p}}
\def\ss{{\bf s}}
\def\tt{{\Theta}}

\def\eb{{\mathbbm{e}}}
\def\rb{{\mathbbm{r}}}

\def\PP{\hbox{\bf P}}

\def\EE{\hbox{\bf E}}

\def\YY{\hbox{\bf Y}}

\def\OO{\mbox{\boldmath$\Omega$}}

\def\llbr{[\hspace{-.10em} [ }
\def\rrbr{ ] \hspace{-.10em}]}

\def\be{\begin{equation}}
\def\ee{\end{equation}}
\def\ba{\begin{eqnarray*}}
\def\ea{\end{eqnarray*}}
\def\ov{\overline}
\def\un{\underline}
\def\wh{\widehat}
\def\wt{\widetilde}

\def\la{\longrightarrow}
\def\lmt{\longmapsto}
\def\da{\downarrow}

\def\deg{{\rm deg}}

\def\proof{\noindent{\textsc{Proof}\hspace{0.1cm}:\hspace{0.15cm}}}

\newtheorem{theorem}{Theorem}[section]
\newtheorem{lemma}[theorem]{Lemma}
\newtheorem{proposition}[theorem]{Proposition}
\newtheorem{corollary}[theorem]{Corollary}

\begin{document}

\title{\bf Asymptotics for rooted planar maps and scaling limits of two-type spatial trees}
\author{Mathilde {\sc Weill}\footnote{DMA-ENS~: 45 rue d'Ulm, 75005 Paris, France -- e-mail~: weill@dma.ens.fr, fax~: (33) 1 44 32 20 80.} \\
{\small \'Ecole normale sup\'erieure de Paris}}
\vspace{2mm}
\date{\today}

\maketitle

\begin{abstract}
We prove some asymptotic results for the radius and the profile of large random bipartite planar maps. Using a bijection due to Bouttier, Di Francesco \& Guitter between rooted bipartite planar maps and certain two-type trees with positive labels, we derive our results from a conditional limit theorem for two-type spatial trees. Finally we apply our estimates to separating vertices of bipartite planar maps~: with probability close to one when $n\to\infty$, a random $2\ka$-angulation with $n$ faces has a separating vertex whose removal disconnects the map into two components each with size greater that $n^{1/2-\vep}$.
\end{abstract}

\section{Introduction}

The main goal of the present work is to investigate asymptotic properties of large rooted bipartite planar maps under the so-called Boltzmann distributions. This setting includes as a special case the asymptotics as $n\to\infty$ of the uniform distribution over rooted $2\ka$-angulations with $n$ faces, for any fixed integer $\ka\geq2$. Boltzmann distributions over planar maps that are both rooted and pointed have been considered recently by Marckert \& Miermont \cite{MaMi} who discuss in particular the profile of distances from the distinguished point in the map. Here we deal with rooted maps and we investigate distances from the root vertex, so that our results do not follow from the ones in \cite{MaMi}, although many statements look similar. The specific results that are obtained in the present work have found applications in the paper \cite{LGcartes}, which investigates scaling limits of large planar maps. 

Let us briefly discuss Boltzmann distributions over rooted bipartite planar maps. We consider a sequence $\q=(q_i)_{i\geq1}$ of weights (nonnegative real numbers) satisfying certain regularity properties. Then, for each integer $n\geq2$, we choose a random rooted bipartite map $M_n$ with $n$ faces whose distribution is specified as follows~: the probability that $M_n$ is equal to a given bipartite planar map $\mm$ is proportional to
$$\prod_{i=1}^nq_{\deg(f_i)/2}$$
where $f_1,\ldots,f_n$ are the faces of $\mm$ and deg($f_i$) is the degree (that is the number of adjacent edges) of the face $f_i$. In particular we may take $q_\ka=1$ and $q_i=0$ for $i\neq\ka$, and we get the uniform distribution over rooted $2\ka$-angulations with $n$ faces.

Theorem \ref{thcartes} below provides asymptotics for the radius and the profile of distances from the root vertex in the random map $M_n$ when $n\to\infty$. The limiting distributions are described in terms of the one-dimensional Brownian snake driven by a normalized excursion. In particular if $R_n$ denotes the radius of $M_n$ (that is the maximal distance from the root), then $n^{-1/4}R_n$ converges to a multiple of the range of the Brownian snake. In the special case of quadrangulations ($q_2=1$ and $q_i=0$ for $i\neq2$), these results were obtained earlier by Chassaing \& Schaeffer \cite{ChS}. As was mentioned above, very similar results have been obtained by Marckert \& Miermont \cite{MaMi} in the setting of Boltzmann distributions over rooted pointed bipartite planar maps, but considering distances from the distinguished point rather than from the root.
 
Similarly as in \cite{ChS} or \cite{MaMi}, bijection between trees and maps serve as a major tool in our approach. In the case of quadrangulations, these bijections were studied by Cori \& Vauquelin \cite{CV} and then by Schaeffer \cite{S}. They have been recently extended to bipartite planar maps by Bouttier, di Francesco \& Guitter \cite{BdFG}. More precisely, Bouttier, di Francesco \& Guitter show that bipartite planar maps are in one-to-one correspondence with well-labelled mobiles, where a well-labelled mobile is a two-type spatial tree whose vertices are assigned positive labels satisfying certain compatibility conditions (see section \ref{secbij} for a precise definition). This bijection has the nice feature that labels in the mobile correspond to distances from the root in the map. Then the above mentioned asymptotics for random maps reduce to a limit theorem for well-labelled mobiles, which is stated as Theorem \ref{thconvcond} below. This statement can be viewed as a conditional version of Theorem 11 in \cite{MaMi}. The fact that \cite{MaMi} deals with distances from the distinguished point in the map (rather than from the root) makes it possible there to drop the positivity constraint on labels. In the present work this constraint makes the proof significantly more difficult. We rely on some ideas from Le Gall \cite{LG} who established a similar conditional theorem for well-labelled trees. Although many arguments in Section \ref{secthconvcond} below are analogous to the ones in \cite{LG}, there are significant additional difficulties because we deal with two-type trees and we condition on the number of vertices of type $1$ rather than on the total number of vertices.

A key step in the proof of Theorem \ref{thconvcond} consists in the derivation of estimates for the probability that a two-type spatial tree remains on the positive half-line. As another application of these estimates, we derive some information about separating vertices of uniform $2\ka$-angulations. We show that with a probability close to one when $n\to\infty$ a random rooted $2\ka$-angulation with $n$ faces will have a vertex whose removal disconnects the map into two components both having size greater that $n^{1/2-\vep}$. Related combinatorial results are obtained in \cite{Airy}. More precisely, in a variety of different models, Proposition 5 in \cite{Airy} asserts that the second largest nonseparable component of a random map of size $n$ has size at most $O(n^{2/3})$. This suggests that $n^{1/2-\vep}$ in our result could be replaced by $n^{2/3-\vep}$.

The paper is organized as follows. In section \ref{secprelim}, we recall some preliminary results and we state our asymptotics for large random rooted planar maps. Section \ref{secthconvcond} is devoted to the proof of Theorem \ref{thconvcond} and to the derivation of Theorem \ref{thcartes} from asymptotics for well-labelled mobiles. Finally Section \ref{secdisc} discusses the application to separating vertices of uniform $2\ka$-angulations. 
 
\section{Preliminaries}\label{secprelim}

\subsection{Boltzmann laws on planar maps}\label{defmaps}

A {\em planar map} is a proper embedding, without edge crossings, of a connected graph in the 2-dimensional sphere $\S^2$. Loops and multiple edges are allowed. A planar map is said to be {\em bipartite} if all its faces have even degree. In this paper, we will only be concerned with bipartite maps. The set of vertices will always be equipped with the graph distance~: if $a$ and $a'$ are two vertices, $d(a,a')$ is the minimal number of edges on a path from $a$ to $a'$. If $M$ is a planar map, we write $\f_M$ for the set of its faces, and $\v_M$ for the set of its vertices.

A {\em pointed} planar map is a pair $(M,\tau)$ where $M$ is a planar map and $\tau$ is a distinguished vertex. Note that, since $M$ is bipartite, if $a$ and $a'$ are two neighbouring vertices, then we have $|d(\tau,a)-d(\tau,a')|=1$. A {\em rooted} planar map is a pair $(M,\vec{e}\,)$ where $M$ is a planar map and $\vec{e}$ is a distinguished oriented edge. The origin of $\vec{e}$ is called the root vertex. At last, a {\em rooted pointed} planar map is a triple $(M,e,\tau)$ where $(M,\tau)$ is a pointed planar map and $e$ is a distinguished non-oriented edge. We can always orient $e$ in such a way that its origin $a$ and its end point $a'$ satisfy $d(\tau,a')=d(\tau,a)+1$. Note that a rooted planar map can be interpreted as a rooted pointed planar map by choosing the root vertex as the distinguished point.

Two pointed maps (resp.~two rooted maps, two rooted pointed maps) are identified if there exists an orientation-preserving homeomorphism of the sphere that sends the first map to the second one and preserves the distinguished point (resp.~the root edge, the distinguished point and the root edge). Let us denote by $\m_p$ (resp.~$\m_r$, $\m_{r,p}$) the set of all pointed bipartite maps (resp.~the set of all rooted bipartite maps, the set of all rooted pointed bipartite maps) up to the preceding indentification.

Let us recall some definitions and propositions that can be found in \cite{MaMi}. Let $\q=(q_i,i\geq1)$ be a sequence of nonnegative weights such that $q_i>0$ for at least one $i>1$. For any planar map $M$, we define $W_\q(M)$ by
$$W_\q(M)=\prod_{f\in\f_M}q_{\deg(f)/2},$$
where we have written ${\rm deg}(f)$ for the degree of the face $f$. We require $\q$ to be admissible that is
$$Z_\q=\sum_{M\in\m_{r,p}}W_\q(M)<\infty.$$
Note that the sum is over the set $\m_{r,p}$ of all rooted pointed bipartite planar maps which is countable thanks to the identification that was explained above. For $k\geq1$, we set $N(k)={2k-1 \choose k-1}$. For every weight sequence $\q$, we define 
$$f_\q(x)=\sum_{k\geq0}N(k+1)q_{k+1}x^k,\;x\geq0.$$
Let $R_\q$ be the radius of convergence of the power series $f_\q$. Consider the equation 
\be\label{eqq}
f_\q(x)=1-x^{-1},\;\,x>0.
\ee
From Proposition 1 in \cite{MaMi}, a sequence $\q$ is admissible if and only if equation (\ref{eqq}) has at least one solution, and then $Z_\q$ is the solution of (\ref{eqq}) that satisfies $Z_\q^2f'_\q(Z_\q)\leq1$. An admissible weight sequence $\q$ is said to be {\em critical} if it satisfies 
$$\left(Z_\q\right)^2f'_\q\left(Z_\q\right)=1,$$
which means that the graphs of the functions $x\lmt f_\q(x)$ and $x\lmt 1-1/x$ are tangent at the left of $x=Z_\q$. Furthermore, if $Z_\q<R_\q$, then $\q$ is said to be {\em regular critical}. This means that the graphs are tangent both at the left and at the right of $Z_\q$. In what follows, we will only be concerned with regular critical weight sequences.

Let $\q$ be a regular critical weight sequence. We define the Boltzmann distribution $\B_\q^{r,p}$ on the set $\m_{r,p}$ by
$$\B_\q^{r,p}(M)=\frac{W_\q(M)}{Z_\q}.$$
Let us now define $Z_\q^{(r)}$ by 
$$Z_\q^{(r)}=\sum_{M\in\m_{r}}\prod_{f\in\f_M}q_{\deg(f)/2}.$$
Note that the sum is over the set $\m_r$ of all rooted bipartite planar maps. From the fact that $Z_\q<\infty$ it easily follows that $Z_\q^{(r)}<\infty$. We then define the Boltzmann distribution $\B_\q^{r}$ on the set $\m_{r}$ by
$$\B_\q^{r}(M)=\frac{W_\q(M)}{Z_\q^{(r)}}.$$

Let us turn to the special case of {\em $2\ka$-angulations}. A $2\ka$-angulation is a bipartite planar map such that all faces have a degree equal to $2\kappa$. If $\kappa=2$, we recognize the well-known {\em quadrangulations}. Let us set
$$\al_\ka=\frac{(\ka-1)^{\ka-1}}{\ka^\ka N(\ka)}.$$ 
We denote by $\q_\kappa$ the weight sequence defined by $q_\kappa=\al_\ka$ and $q_i=0$ for every $i\in\N\setminus\{\kappa\}$. It is proved in Section 1.5 of \cite{MaMi} that $\q_\ka$ is a regular critical weight sequence, and
$$Z_{\q_\ka}=\frac{\ka}{\ka-1}.$$
For every $n\geq1$, we denote by $\U^n_\ka$ (resp.~$\ov{\U}^{\,n}_\ka$) the uniform distribution on the set of all rooted pointed $2\ka$-angulations with $n$ faces (resp.~on the set of all rooted $2\ka$-angulations with $n$ faces). We have 
\ba
\B^{r,p}_{\q_\kappa}(\cdot\mid\#\f_M=n)&=&\U^n_\ka,\\
\B^{r}_{\q_\kappa}(\cdot\mid\#\f_M=n)&=&\ov{\U}^{\,n}_\ka.
\ea

\subsection{Two-type spatial Galton-Watson trees}\label{twotypetrees}

We start with some formalism for discrete trees. Set
$$\u=\bigcup_{n\geq0}\N^n,$$ 
where $\N=\{1,2,3,\ldots\}$ and by convention $\N^0=\{\emptyset\}$. An element of $\,\u$ is a sequence $u=u^1\ldots u^n$, and we set $|u|=n$ so that $|u|$ represents the generation of $u$. In particular, $\vert\emptyset\vert=0$. If $u=u^1\ldots u^n$ and $v=v^1\ldots v^m$ belong to $\,\u$, we write $uv=u^1\ldots u^nv^1\ldots v^m$ for the concatenation of $u$ and $v$. In particular, $\emptyset u=u\emptyset=u$. If $v$ is of the form $v=uj$ for $u\in\u$ and $j\in\N$, we say that $v$ is a child of $u$, or that $u$ is the father of $v$, and we write $u=\check{v}$. More generally if $v$ is of the form $v=uw$ for $u,w\in\u$, we say that $v$ is a descendant of $u$, or that $u$ is an ancestor of $v$. The set $\u$ comes with the natural lexicographical order such that $u\preccurlyeq v$ if either $u$ is an ancestor of $v$, or if $u=wa$ and $v=wb$ with $a\in\u^\ast$ and $b\in\u^\ast$ satisfying $a^1<b^1$, where we have set $\u^\ast=\u\setminus\{\emptyset\}$. And we write $u\prec v$ if $u\preccurlyeq v$ and $u\neq v$.

A plane tree $\t$ is a finite subset of $\,\u$ such that 
\begin{enumerate}
\item[(i)] $\emptyset\in\t$,
\item[(ii)] $u\in\t\setminus\{\emptyset\}\Rightarrow\check{u}\in\t$,
\item[(iii)] for every $u\in\t$, there exists a number $k_u(\t)\geq0$ such that $uj\in\t$ if and only if $1\leq j\leq k_u(\t)$.
\end{enumerate}
We denote by $\T$ the set of all plane trees. 

Let $\t$ be a plane tree and let $\zeta=\#\t-1$. The {\em search-depth sequence} of $\t$ is the sequence $u_0,u_1,\ldots,u_{2\zeta}$ of vertices of $\t$ wich is obtained by induction as follows. First $u_0=\emptyset$, and then for every $i\in\{0,1,\ldots,2\zeta-1\}$, $u_{i+1}$ is either the first child of $u_i$ that has not yet appeared in the sequence $u_0,u_1,\ldots,u_i$, or the father of $u_i$ if all children of $u_i$ already appear in the sequence $u_0,u_1,\ldots,u_i$. It is easy to verify that $u_{2\zeta}=\emptyset$ and that all vertices of $\t$ appear in the sequence $u_0,u_1,\ldots,u_{2\zeta}$ (of course some of them appear more that once). We can now define the {\em contour function} of $\t$. For every $k\in\{0,1,\ldots,2\zeta\}$, we let $C(k)$ denote the distance from the root of the vertex $u_{k}$. We extend the definition of $C$ to the line interval $[0,2\zeta]$ by interpolating linearly between successive integers. Clearly $\t$ is uniquely determined by its contour function $C$.

A discrete spatial tree is a pair $(\t,U)$ where $\t\in\T$ and $U=(U_v,v\in\t)$ is a mapping from the set $\t$ into $\R$. If $v$ is a vertex of $\t$, we say that $U_v$ is the {\em label} of $v$. We denote by $\Omega$ the set of all discrete spatial trees. If $(\t,U)\in\Omega$ we define the {\em spatial contour function} of $(\t,U)$ as follows. First if $k$ is an integer, we put $V(k)=U_{u_k}$ with the preceding notation. We then complete the definition of $V$ by interpolating linearly between successive integers. Clearly $(\t,U)$ is uniquely determined by the pair $(C,V)$.

Let $(\t,U)\in\Omega$. We interpret $(\t,U)$ as a two-type (spatial) tree by declaring that vertices of even generations are of type $0$ and vertices of odd generations are of type $1$. We then set 
\ba
\t^0&=&\{u\in\t:|u|\;{\rm is}\;{\rm even}\},\\
\t^1&=&\{u\in\t:|u|\;{\rm is}\;{\rm odd}\}.
\ea 

Let us turn to random trees. We want to consider a particular family of two-type Galton-Watson trees, in which vertices of type $0$ only give birth to vertices of type $1$ and vice-versa. Let $\mu=(\mu_0,\mu_1)$ be a pair of offspring distributions, that is a pair of probability distributions on $\Z_+$. If $m_0$ and $m_1$ are the respective means of $\mu_0$ and $\mu_1$ we assume that $m_0m_1\leq1$ and we exclude the trivial case $\mu_0=\mu_1=\dl_1$ where $\dl_1$ stands for the Dirac mass at $1$. We denote by $P_\mu$ the law of a two-type Galton-Watson tree with offspring distribution $\mu$, meaning that for every $\tr\in\T$,
$$P_\mu(\tr)=\prod_{u\in\tr^0}\mu_0\left(k_u(\tr)\right)\prod_{u\in\tr^1}\mu_1\left(k_u(\tr)\right),$$
where $\tr^0$ (resp.~$\tr^1$) is as above the set of all vertices of $\tr$ with even (resp.~odd) generation. The fact that this formula defines a probability measure on $\T$ is justified in \cite{MaMi}.

Let us now recall from \cite{MaMi} how one can couple plane trees with a spatial displacement in order to turn them into random elements of $\Omega$. To this end, let $\nu_0^k,\nu_1^k$ be probability distributions on $\R^k$ for every $k\geq1$. We set $\nu=(\nu_0^k,\nu_1^k)_{k\geq1}$.  For every $\t\in\T$ and $x\in\R$, we denote by $R_{\nu,x}(\t,\d U)$ the probability measure on $\R^\t$ which is characterized as follows. Let $(\YY_u,u\in\t)$ be a family of independent random variables such that for $u\in\t$ with $k_u(\t)=k$, $\YY_u=(Y_{u1},\ldots,Y_{uk})$ is distributed according to $\nu_0^k$ if $u\in\t^0$ and according to $\nu_1^k$ if $u\in\t^1$. We set $X_\emptyset=x$ and for every $v\in\t\setminus\{\emptyset\}$,
$$X_v=x+\sum_{u\in]\emptyset,v]}Y_u,$$
where $]\emptyset,v]$ is the set of all ancestors of $v$ distinct from the root $\emptyset$. Then $R_{\nu,x}(\t,\d U)$ is the law of $(X_v,v\in\t)$. We finally define for every $x\in\R$ a probability measure $\P_{\mu,\nu,x}$ on $\Omega$ by setting
$$\P_{\mu,\nu,x}(\d\t\,\d U)=P_\mu(\d\t)R_{\nu,x}(\t,\d U).$$

\subsection{The Brownian snake and the conditioned Brownian snake}\label{serpent}

Let $x\in\R$. The Brownian snake with initial point $x$ is a pair $(\e,\r^x)$, where $\e=(\e(s),0\leq s\leq1)$ is a normalized Brownian excursion and $\r^x=(\r^x(s),0\leq s\leq 1)$ is a real-valued process such that, conditionally given $\e$, $\r^x$ is Gaussian with mean and covariance given by
\begin{enumerate}
\item[{$\bullet$}] $\EE[\r^x(s)]=x$ for every $s\in[0,1]$,
\item[{$\bullet$}] ${\bf Cov}(\r^x(s),\r^x(s'))=\displaystyle{\inf_{s\leq 
t\leq s'}\e(t)}$ for every $0\leq s\leq s'\leq 1$.
\end{enumerate}
We know from \cite{Zu} that $\r^x$ admits a continuous modification. From now on we consider only this modification. In the terminology of \cite{Zu} $\r^x$ is the terminal point process of the one-dimensional Brownian snake driven by the normalized Brownian excursion $\e$ and with initial point $x$. 

Write $\PP$ for the probability measure under which the collection $(\e,\r^x)_{x\in\R}$ is defined. Note that for every $x>0$, we have
$$\PP\left(\inf_{s\in[0,1]}\r^x(s)\geq0\right)>0.$$
We may then define for every $x>0$ a pair $(\ov{\e}^x,\ov{\r}^x)$ which is distributed as the pair $(\e,\r^x)$ under the conditioning that $\inf_{s\in[0,1]}\r^x(s)\geq0$. 

We equip $C([0,1],\R)^2$ with the norm $\|(f,g)\|=\|f\|_u\vee\|g\|_u$ where $\|f\|_u$ stands for the supremum norm of $f$. The following theorem is Theorem Theorem 1.1 in \cite{LGW}.

\begin{theorem}
There exists a pair $(\ov{\e}^0,\ov{\r}^0)$ such that $(\ov{\e}^x,\ov{\r}^x)$ converges in distribution as $x\da0$ towards $(\ov{\e}^0,\ov{\r}^0)$.
\end{theorem} 
The pair $(\ov{\e}^0,\ov{\r}^0)$ is the so-called conditioned Brownian snake with initial point $0$.

Theorem 1.2 in \cite{LGW} provides a useful construction of the conditioned object $(\ov{\e}^0,\ov{\r}^0)$ from the unconditioned one $(\e,\r^0)$. In order to present this construction, first recall that there is a.s.~a unique $s_*$ in $(0,1)$ such that
$$\r^0(s_*)=\inf_{s\in[0,1]}\r^0(s)$$
(see Lemma 16 in \cite{MaMo2} or Proposition 2.5 in \cite{LGW}). For every $s\in[0,\infty)$, write $\{s\}$ for the fractional part of $s$. According to Theorem 1.2 in \cite{LGW}, the conditioned snake $(\ov{\e}^0,\ov{\r}^0)$ may be constructed explicitly as follows~: for every $s\in[0,1]$,
\ba
\ov{\e}^0(s)&=&\e({s_*})+\e(\{s_*+s\})-2\,\inf_{s\wedge\{s_*+s\}\leq t\leq s\vee\{s_*+s\}}\e(t),\\
\ov{\r}^0(s)&=&\r^0(\{s_*+s\})-\r^0({s_*}).
\ea

\subsection{The Bouttier-di Francesco-Guitter bijection}\label{secbij}

We start with a definition. A (rooted) {\em mobile} is a two-type spatial tree $(\t,U)$ whose labels $U_v$ only take integer values and such that the following properties hold~:
\begin{enumerate}
\item[(a)] $U_v=U_{\check{v}}$ for every $v\in\t^1$.
\item[(b)] Let $v\in\t^1$ such that $k=k_v(\t)\geq1$. Let $v_{(0)}=\check{v}$ be the father of $v$ and let $v_{(j)}=vj$ for every $j\in\{1,\ldots,k\}$. Then for every $j\in\{0,\ldots,k\}$, 
$$U_{v_{(j+1)}}\geq U_{v_{(j)}}-1,$$ 
where by convention $v_{(k+1)}=v_{(0)}$.
\end{enumerate} 
Furthermore, if $U_v\geq1$ for every $v\in\t$, then we say that $(\t,U)$ is a {\em well-labelled} mobile. 

Let $\T^{\rm mob}_1$ denotes the set of all mobiles such that $U_\emptyset=1$. We will now describe the Bouttier-di Francesco-Guitter bijection from $\T^{\rm mob}_1$ onto $\m_{r,p}$. This bijection can be found in section 2 in \cite{BdFG}. Note that \cite{BdFG} deals with pointed planar maps rather than with rooted pointed planar maps. It is however easy to verify that the results described below are simple consequences of \cite{BdFG}. 

Let $(\t,U)\in\T^{\rm mob}_1$. Recall that $\zeta=\#\t-1$. Let $u_0,u_1,\ldots,u_{2\zeta}$ be the search-depth sequence of $\t$. It is immediate to see that $u_k\in\t^0$ if $k$ is even and that $u_k\in\t^1$ if $k$ is odd. The search-depth sequence of $\t^0$ is the sequence $w_0,w_1,\ldots,w_\zeta$ defined by $w_k=u_{2k}$ for every $k\in\{0,1,\ldots,\zeta\}$. Notice that $w_0=w_\zeta=\emptyset$. Although $(\t,U)$ is not necesseraly well labelled, we may set for every $v\in\t$, 
$$U^+_v=U_v-\min\{U_w:w\in\t\}+1,$$ 
and then $(\t,U^+)$ is a well-labelled mobile. Notice that $\min\{U^+_v:v\in\t\}=1$.

Suppose that the tree $\t$ is drawn in the plane and add an extra vertex $\dd$. We associate with $(\t,U^+)$ a bipartite planar map whose set of vertices is 
$$\t^0\cup\{\dd\},$$
and whose edges are obtained by the following device~: for every $k\in\{0,1,\ldots,\zeta\}$,
\begin{enumerate}
\item[$\bullet$] if $U^+_{w_k}=1$, draw an edge between $w_k$ and $\dd$~;
\item[$\bullet$] if $U^+_{w_k}\geq2$, draw an edge between $w_k$ and the first vertex in the sequence $w_{k+1},\ldots,w_{\zeta-1}$, $w_0,w_1,\ldots,w_{k-1}$ whose label is $U^+_{w_k}-1$.
\end{enumerate}
Notice that condition (b) in the definition of a mobile entails that $U^+_{w_{k+1}}\geq U^+_{w_k}-1$ for every $k\in\{0,1,\ldots,\zeta-1\}$ and recall that $\min\{U^+_{w_0},U^+_{w_1},\ldots,U^+_{w_{\zeta-1}}\}=1$. The preceding properties ensure that whenever $U^+_{w_k}\geq2$ there is at least one vertex among $w_{k+1},\ldots,w_{\zeta-1},w_0,\ldots,w_{k-1}$ with label $U^+_{w_{k}}-1$. The construction can be made in such a way that edges do not intersect (see section 2 in \cite{BdFG} for an example). The resulting planar graph is a bipartite planar map. We view this map as a rooted pointed planar map by declaring that the distinguished vertex is $\dd$ and that the root edge is the one corresponding to $k=0$ in the preceding construction.

It follows from \cite{BdFG} that the preceding construction yields a bijection $\Psi_{r,p}$ between $\T^{\rm mob}_1$ and $\m_{r,p}$. Furthermore it is not difficult to see that $\Psi_{r,p}$ satisfies the following two properties~: let $(\t,U)\in\T^{\rm mob}_1$ and let $M=\Psi_{r,p}((\t,U))$,
\begin{enumerate}
\item[(i)] for every $k\geq1$, the set $\{f\in\f_M:\deg(f)=2k\}$ is in one-to-one correspondence with the set $\{v\in\t^1:k_v(\t)=k-1\}$,
\item[(ii)] for every $l\geq1$, the set $\{a\in\v_M:d(\dd,a)=l\}$ is in one-to-one correspondence with the set $\{v\in\t^0:U_v-\min\{U_w:w\in\t\}+1=l\}$.
\end{enumerate} 

We observe that if $(\t,U)$ is a well-labelled mobile then $U^+_v=U_v$ for every $v\in\t$. In particular $U^+_\emptyset=1$. This implies that the root edge of the planar map $\Psi_{r,p}((\t,U))$ contains the distinguished point $\dd$. Then $\Psi_{r,p}((\t,U))$ can be identified to a rooted planar map, whose root is an oriented edge between the root vertex $\dd$ and $w_0$. Write $\ov{\T}^{\rm mob}_1$ for the set of all well-labelled mobiles such that $U_\emptyset=1$. Thus $\Psi_{r,p}$ induces a bijection $\Psi_r$ from the set $\ov{\T}^{\rm mob}_1$ onto the set $\m_r$. Furthermore $\Psi_r$ satisfies the following two properties~: let $(\t,U)\in\ov{\T}^{\rm mob}_1$ and let $M=\Psi_{r}((\t,U))$,
\begin{enumerate}
\item[(i)] for every $k\geq1$, the set $\{f\in\f_M:\deg(f)=2k\}$ is in one-to-one correspondence with the set $\{v\in\t^1:k_v(\t)=k-1\}$,
\item[(ii)] for every $l\geq1$, the set $\{a\in\v_M:d(\dd,a)=l\}$ is in one-to-one correspondence with the set $\{v\in\t^0:U_v=l\}$.
\end{enumerate} 

\subsection{Boltzmann distribution on two-type spatial trees} 

Let $\q$ be a regular critical weight sequence. We recall the following definitions from \cite{MaMi}. Let $\mu_0^\q$ be the geometric distribution with parameter $f_\q(Z_\q)$ that is
$$\mu_0^\q(k)=Z_\q^{-1}f_\q(Z_\q)^k,\;\,k\geq0,$$
and let $\mu_1^\q$ be the probability measure defined by
$$\mu_1^\q(k)=\frac{Z_\q^kN(k+1)q_{k+1}}{f_\q(Z_\q)},\;\,k\geq0.$$ 
From \cite{MaMi}, we know that $\mu_1$ has small exponential moments, and that the two-type Galton-Watson tree associated with $\mu^\q=(\mu_0^\q,\mu_1^\q)$ is critical.  

Also, for every $k\geq0$, let $\nu_0^k$ be the Dirac mass at $0\in\R^k$ and let $\nu_1^k$ be the uniform distribution on the set $A_k$ defined by
$$A_k=\left\{(x_1,\ldots,x_k)\in\Z^k:x_1\geq-1,x_2-x_1\geq-1,\ldots,x_k-x_{k-1}\geq-1,-x_k\geq-1\right\}.$$ 
We can say equivalently that $\nu_1^k$ is the law of $(X_1,\ldots,X_1+\ldots+X_k)$ where $(X_1,\ldots,X_{k+1})$ is uniformly distributed on the set $B_k$ defined by
$$B_k=\left\{(x_1,\ldots,x_{k+1})\in\{-1,0,1,2,\ldots\}^{k+1}:x_1+\ldots+x_{k+1}=0\right\}.$$ 
Notice that $\#A_k=\#B_k=N(k+1)$. We set $\nu=((\nu_0^k,\nu_1^k))_{k\geq1}$. The following result is Proposition 10 in \cite{MaMi}. However, we provide a short proof for the sake of completeness.

\begin{proposition}\label{loiBolGW}
The Boltzmann distribution $\B^{r,p}_\q$ is the image of the probability measure $\P_{\mu^\q,\nu,1}$ under the mapping $\Psi_{r,p}$. 
\end{proposition}

\proof By construction, the probability measure $\P_{\mu^\q,\nu,1}$ is supported on the set $\T^{\rm mob}_1$. Let $(\tt,\uu)\in\T^{\rm mob}_1$. We have by the choice of $\nu$,
\begin{eqnarray*}
\P_{\mu^\q,\nu,1}\left((\tr,\uu)\right)&=&P_{\mu^\q}(\tr)R_{\nu,1}(\tr,\{\uu\})\\
&=&\Big(\prod_{v\in\tr^{1}}N(k_v(\tr)+1)\Big)^{-1}P_{\mu^\q}(\tr).
\end{eqnarray*}
Now,
\begin{eqnarray*} 
P_{\mu^\q}(\tr)&=&\prod_{v\in\tr^0}\mu^\q_0(k_v(\tr))\prod_{v\in\tr^1}\mu_1^\q(k_v(\tr))\\
&=&\prod_{v\in\tr^0}\left(Z_\q^{-1}f_\q(Z_\q)^{k_v(\tr)}\right)\prod_{v\in\tr^{1}}\frac{Z_\q^{k_v(\tr)}N(k_v(\tr)+1)q_{k_v(\tr)+1}}{f_\q(Z_\q)}\\
&=&Z_\q^{-\#\tr^{0}}f_\q(Z_\q)^{\#\tr^{1}}Z_\q^{\#\tr^{0}-1}f_\q(Z_\q)^{-\#\tr^{1}}\prod_{v\in\tr^{1}}N(k_v(\tr)+1)\prod_{v\in\tr^{1}}q_{k_v(\tr)+1}\\
&=&Z_\q^{-1}\prod_{v\in\tr^{1}}q_{k_v(\tr)+1}\prod_{v\in\tr^{1}}N(k_v(\tr)+1),
\end{eqnarray*}
so that we arrive at
$$\P_{\mu^\q,\nu,1}\left((\tr,\uu)\right)=Z_\q^{-1}\prod_{v\in\tr^{1}}q_{k_v(\tr)+1}.$$
We set $\mm=\Psi_{r,p}((\tr,\uu))$. We have from the property (ii) satisfied by $\Psi_{r,p}$,
$$Z_\q^{-1}\prod_{v\in\tr^{1}}q_{k_v(\tr)+1}=\B^{r,p}_\q(\mm),$$ 
which leads us to the desired result. \cq

Let us introduce some notation. As $\mu_0^\q(1)>0$, we have $P_\mu(\t^1=n)>0$ for every $n\geq1$. Then we may define, for every $n\geq1$ and $x\in\R$, 
\ba
P_{\mu^\q}^n&=&P_{\mu^\q}\left(\cdot\mid\#\t^1=n\right),\\
\P^n_{\mu^\q,\nu,x}&=&\P_{\mu^\q,\nu,x}\left(\cdot\mid\#\t^1=n\right).
\ea
Furthermore, we set for every $(\t,U)\in\Omega$,
$$\un{U}=\min\left\{U_v:v\in\t^0\setminus\{\emptyset\}\right\},$$ 
with the convention $\min\emptyset=\infty$. Finally we define for every $n\geq1$ and $x\geq0$, 
\ba
\ov{\P}_{\mu^\q,\nu,x}&=&\P_{\mu^\q,\nu,x}(\cdot\mid\un{U}>0),\\
\ov{\P}^{\,n}_{\mu^\q,\nu,x}&=&\ov{\P}_{\mu^\q,\nu,x}\left(\cdot\mid\#\t^1=n\right).
\ea

\begin{corollary}\label{loiimage}
The probability measure $\B^{r,p}_\q(\cdot\mid\#\f_M=n)$ is the image of $\P^n_{\mu^\q,\nu,1}$ under the mapping $\Psi_{r,p}$. The probability measure $\B^{r}_\q$ is the image of $\ov{\P}_{\mu^\q,\nu,1}$ under the mapping $\Psi_r$. The probability measure $\B^{r}_\q(\cdot\mid\#\f_M=n)$ is the image of $\ov{\P}^{\,n}_{\mu^\q,\nu,1}$ under the mapping $\Psi_r$. 
\end{corollary}

\proof The first assertion is a simple consequence of Proposition \ref{loiBolGW} together with the property (i) satisfied by $\Psi_{r,p}$. Recall from section \ref{defmaps} that we can identify the set $\m_r$ to a subset of $\m_{r,p}$ in the following way. Let $\Upsilon:\m_r\la\m_{r,p}$ be the mapping defined by $\Upsilon((M,\vec{e}\,))=(M,e,o)$ for every $(M,\vec{e}\,)\in\m_r$, where $o$ denotes the root vertex of the map $(M,\vec{e}\,)$. We easily check that $\B^{r,p}_\q(\cdot\mid M\in\Upsilon(\m_r))$ is the image of $\B^r_\q$ under the mapping $\Upsilon$. This together with Proposition \ref{loiBolGW} yields the second assertion. The third assertion follows.  \cq

At last, if $\q=\q_\ka$, we set $\mu_0^\ka=\mu_0^{\q_\ka}$, $\mu_1^\ka=\mu_1^{\q_\ka}$ and $\mu^\ka=(\mu_0^\ka,\mu_1^\ka)$. We then verify that $\mu_0^\ka$ is the geometric distribution with parameter $1/\ka$ and that $\mu_1^\ka$ is the Dirac mass at $\ka-1$. Recall the notation $\U^n_\ka$ and $\ov{\U}^{\,n}_\ka$. 
\begin{corollary}\label{loiimageka}
The probability measure $\U^n_\ka$ is the image of $\P^n_{\mu^\ka,\nu,1}$ under the mapping $\Psi_{r,p}$. The probability measure $\ov{\U}^{\,n}_\ka$ is the image of $\ov{\P}^{\,n}_{\mu^\ka,\nu,1}$ under the mapping $\Psi_r$. 
\end{corollary}

\subsection{Statement of the main result}

We first need to introduce some notation. Let $M\in\m_r$. We denote by $o$ its root vertex. The radius $\rad_M$ is the maximal distance between $o$ and another vertex of $M$ that is
$$\rad_M=\max\{d(o,a):a\in\v_M\}.$$
The normalized profile of $M$ is the probability measure $\lm_M$ on $\{0,1,2,\ldots\}$  defined by
$$\lm_M(k)=\frac{\#\{a\in\v_M:d(o,a)=k\}}{\#\v_M},\;k\geq0.$$
Note that $\rad_M$ is the supremum of the support of $\lm_M$. It is also convenient to introduce the rescaled profile. If $M$ has $n$ faces, this is the probability measure on $\R_+$ defined by
$$\lm_M^{(n)}(A)=\lm_M\left(n^{1/4}A\right)$$
for any Borel subset $A$ of $\R_+$. At last, if $\q$ is a regular critical weight sequence, we set
$$\rho_\q=2+Z_\q^3f''_\q(Z_\q).$$

Recall from section \ref{serpent} that $(\e,\r^0)$ denotes the Brownian snake with initial point $0$.

\begin{theorem}\label{thcartes}
Let $\q$ be a regular critical weight sequence. 
\begin{enumerate}
\item[(i)] The law of $n^{-1/4}\,\rad_M$ under the probability measure $\B^r_\q(\cdot\mid\#\f_M=n)$ converges as $n\to\infty$ to the law of the random variable
$$\left(\frac{4\rho_\q}{9(Z_\q-1)}\right)^{1/4}\left(\sup_{0\leq s\leq1}\r^0(s)-\inf_{0\leq s\leq1}\r^0(s)\right).$$
\item[(ii)] The law of the random measure $\lm_M^{(n)}$ under the probability measure $\B^r_\q(\cdot\mid\#\f_M=n)$ converges as $n\to\infty$ to the law of the random probability measure $\I$ defined by
$$\langle\I,g\rangle=\int_0^1g\left(\left(\frac{4\rho_\q}{9(Z_\q-1)}\right)^{1/4}\left(\r^0(t)-\inf_{0\leq s\leq1}\r^0(s)\right)\right)\d t.$$
\item[(iii)] The law of the rescaled distance $n^{-1/4}\,d(o,a)$ where $a$ is a vertex chosen uniformly at random among all vertices of $M$, under the probability measure $\B^r_\q(\cdot\mid\#\f_M=n)$ converges as $n\to\infty$ to the law of the random variable
$$\left(\frac{4\rho_\q}{9(Z_\q-1)}\right)^{1/4}\left(\sup_{0\leq s\leq1}\r^0(s)\right).$$
\end{enumerate}
\end{theorem}
In the case $\q=\q_\ka$, the constant appearing in Theorem \ref{thcartes} is $(4\ka(\ka-1)/9)^{1/4}$. It is equal to $(8/9)^{1/4}$ when $\ka=2$. The results stated in Theorem \ref{thcartes} in the special case $\q=\q_2$ were obtained by Chassaing \& Schaeffer \cite{ChS} (see also Theorem 8.2 in \cite{LG}).

Obviously Theorem \ref{thcartes} is related to Theorem 3 proved by Marckert \& Miermont \cite{MaMi}. Note however that \cite{MaMi} deals with rooted pointed maps instead of rooted maps as we do and studies distances from the distinguished point of the map rather than from the root vertex.

\section{A conditional limit theorem for two-type spatial trees}\label{secthconvcond}

Recall first some notation. Let $\q$ be a regular critical weight sequence, let $\mu^\q=(\mu^\q_0,\mu^\q_1)$ be the pair of offspring distributions associated with $\q$ and let $\nu=(\nu^k_0,\nu^k_1)_{k\geq1}$ be the family of probability measures defined before Proposition \ref{loiBolGW}.

If $(\t,U)\in\Omega$, we denote by $C$ its contour function and by $V$ its spatial contour function. 

Recall that $C([0,1],\R)^2$ is equipped with the norm $\|(f,g)\|=\|f\|_u\vee\|g\|_u$. The following result is a special case of Theorem 11 in \cite{MaMi}.

\begin{theorem}\label{thconv} Let $\q$ be a regular critical weight sequence. The law under $\P_{\mu^\q,\nu,0}^n$ of 
$$\left(\left(\frac{\sqrt{\rho_\q(Z_\q-1)}}{4}\;\frac{C(2(\#\t-1)t)}{n^{1/2}}\right)_{0\leq t\leq1},\left(\left(\frac{9(Z_\q-1)}{4\rho_\q}\right)^{1/4}\frac{V(2(\#\t-1)t)}{n^{1/4}}\right)_{0\leq t\leq1}\right)$$
converges as $n\to\infty$ to the law of $(\e,\r^0)$. The convergence holds in the sense of weak convergence of probability measures on $C([0,1],\R)^2$.
\end{theorem}
Note that Theorem 11 in \cite{MaMi} deals with the so-called height-process instead of the contour process. However, we can deduce Theorem \ref{thconv} from \cite{MaMi} by classical arguments (see e.g.~\cite{LGDEA}).

In this section, we will prove a conditional version of Theorem \ref{thconv}. Before stating this result, we establish a corollary of Theorem \ref{thconv}. To this end we set
\ba
Q_{\mu^\q}&=&P_{\mu^\q}(\cdot\mid k_\emptyset(\t)=1),\\
\Q_{\mu^\q,\nu}&=&\P_{\mu^\q,\nu,0}(\cdot\mid k_\emptyset(\t)=1).
\ea
Notice that this conditioning makes sense since $\mu^\q_0(1)>0$. We may also define for every $n\geq1$,
\ba
Q_{\mu^\q}^n&=&Q_{\mu^\q}\left(\cdot\mid\#\t^1=n\right),\\
\Q_{\mu^\q,\nu}^n&=&\Q_{\mu^\q,\nu}\left(\cdot\mid\#\t^1=n\right).
\ea

\begin{corollary}\label{corconvQ}
Let $\q$ be a regular critical weight sequence. The law under $\Q_{\mu^\q,\nu}^n$ of 
$$\left(\left(\frac{\sqrt{\rho_\q(Z_\q-1)}}{4}\;\frac{C(2(\#\t-1)t)}{n^{1/2}}\right)_{0\leq t\leq1},\left(\left(\frac{9(Z_\q-1)}{4\rho_\q}\right)^{1/4}\frac{V(2(\#\t-1)t)}{n^{1/4}}\right)_{0\leq t\leq1}\right)$$
converges as $n\to\infty$ to the law of $(\e,\r^0)$. The convergence holds in the sense of weak convergence of probability measures on $C([0,1],\R)^2$.
\end{corollary}

\proof We first introduce some notation. If $(\t,U)\in\Omega$ and $w_0\in\t$, we define a spatial tree $(\t^{[w_0]},U^{[w_0]})$ by setting 
$$\t^{[w_0]}=\{v:w_0v\in\t\},$$ 
and for every $v\in\t^{[w_0]}$
$$U^{[w_0]}_v=U_{w_0v}-U_{w_0}.$$ 
Denote by $C^{[w_0]}$ the contour function and by $V^{[w_0]}$ the spatial contour function of $(\t^{[w_0]},U^{[w_0]})$.

As a consequence of Theorem 11 in \cite{MaMi}, the law under $\Q_{\mu^\q,\nu}^n$ of
$$\left(\left(\frac{\sqrt{\rho_\q(Z_\q-1)}}{4}\;\frac{C^{[1]}\left(2\left(\#\t^{[1]}-1\right)t\right)}{n^{1/2}}\right)_{0\leq t\leq1}\hspace{-0.25cm},
\left(\left(\frac{9(Z_\q-1)}{4\rho_\q}\right)^{1/4}\frac{V^{[1]}\left(2\left(\#\t^{[1]}-1\right)t\right)}{n^{1/4}}\right)_{0\leq t\leq1}\right)$$
converges as $n\to\infty$ to the law of $(\e,\r^0)$. We then easily get the desired result.\cq  

Recall from section \ref{serpent} that $(\ov{\e}^0,\ov{\r}^0)$ denotes the conditioned Brownian snake with initial point $0$.

\begin{theorem}\label{thconvcond} Let $\q$ be a regular critical weight sequence. For every $x\geq0$, the law under $\ov{\P}_{\mu^\q,\nu,x}^{\,n}$ of 
$$\left(\left(\frac{\sqrt{\rho_\q(Z_\q-1)}}{4}\;\frac{C(2(\#\t-1)t)}{n^{1/2}}\right)_{0\leq t\leq1},\left(\left(\frac{9(Z_\q-1)}{4\rho_\q}\right)^{1/4}\frac{V(2(\#\t-1)t)}{n^{1/4}}\right)_{0\leq t\leq1}\right)$$
converges as $n\to\infty$ to the law of $(\ov{\e}^0,\ov{\r}^0)$. The convergence holds in the sense of weak convergence of probability measures on $C([0,1],\R)^2$.
\end{theorem}
To prove Theorem \ref{thconvcond}, we will follow the lines of the proof of Theorem 2.2 in \cite{LG}. From now on, we set $\mu=\mu^\q$ to simplify notation.

\subsection{Rerooting two-type spatial trees}

If $\t\in\T$, we say that a vertex $v\in\t$ is a leaf of $\t$ if $k_v(\t)=0$ meaning that $v$ has no child. We denote by $\dd\t$ the set of all leaves of $\t$ and we write $\dd_0\t=\dd\t\cap\t^0$ for the set of leaves of $\t$ which are of type $0$.

Let us recall some notation that can be found in section 3 in \cite{LG}. Recall that $\u^\ast=\u\setminus\{\emptyset\}$. If $v_0\in\u^\ast$ and $\t\in\T$ are such that $v_0\in\t$, we define $k=k(v_0,\t)$ and $l=l(v_0,\t)$ in the following way. Write $\zeta=\#\t-1$ and $u_0,u_1,\ldots,u_{2\zeta}$ for the search-depth sequence of $\t$. Then we set
\ba
k&=&\min\{i\in\{0,1,\ldots,2\zeta\}:u_i=v_0\},\\
l&=&\max\{i\in\{0,1,\ldots,2\zeta\}:u_i=v_0\},
\ea
which means that $k$ is the time of the first visit of $v_0$ in the evolution of the contour of $\t$ and that $l$ is the time of the last visit of $v_0$. Note that $l\geq k$ and that $l=k$ if and only if $v_0\in\dd\t$. For every $t\in[0,2\zeta-(l-k)]$, we set
$$\wh{C}^{(v_0)}(t)=C(k)+C(\llbr k-t\rrbr)-2\inf_{s\in[k\wedge\llbr k-t\rrbr,k\vee\llbr k-t\rrbr]}C(s),$$
where $C$ is the contour function of $\t$ and $\llbr k-t\rrbr$ stands for the unique element of $[0,2\zeta)$ such that $\llbr k-t\rrbr-(k-t)=0$ or $2\zeta$. Then there exists a unique plane tree $\wh{\t}^{(v_0)}\in\T$ whose contour function is $\wh{C}^{(v_0)}$. Informally, $\wh{\t}^{(v_0)}$ is obtained from $\t$ by removing all vertices that are descendants of $v_0$ and by re-rooting the resulting tree at $v_0$. Furthermore, if $v_0=u^1\ldots u^n$, then we see that $\wh{v}_0=1u^n\ldots u^2$ belongs to $\wh{\t}^{(v_0)}$. In fact, $\wh{v}_0$ is the vertex of $\wh{\t}^{(v_0)}$ corresponding to the root of the initial tree. At last notice that $k_\emptyset(\wh{\t}^{(v_0)})=1$.

If $\t\in\T$ and $w_0\in\t$, we set 
$$\t^{(w_0)}=\t\setminus\left\{w_0u\in\t:u\in\u^\ast\right\}.$$ 

The following lemma is an analogue of Lemma 3.1 in \cite{LG} for two-type Galton-Watson trees. Note that in what follows, two-type trees will always be re-rooted at a vertex of type $0$. 

Recall the definition of the probability measure $Q_\mu$.

\begin{lemma}\label{reenracinement} Let $v_0\in\u^\ast$ be of the form $v_0=1u^2\ldots u^{2p}$ for some $p\in\N$. Assume that $Q_\mu(v_0\in\t)>0$. Then the law of the re-rooted tree $\wh{\t}^{(v_0)}$ under $Q_\mu(\cdot\mid v_0\in\t)$ coincides with the law of the tree $\t^{(\wh{v}_0)}$ under $Q_\mu(\cdot\mid\wh{v}_0\in\t)$.
\end{lemma}

\proof We first notice that 
$$Q_\mu(v_0\in\t)=\mu_1\left(\left\{u^2,u^2+1,\ldots\right\}\right)\mu_0\left(\left\{u^3,u^3+1,\ldots\right\}\right)\ldots\mu_1\left(\left\{u^{2p},u^{2p}+1,\ldots\right\}\right),$$
so that 
$$Q_\mu(v_0\in\t)=Q_\mu(\wh{v}_0\in\t)>0.$$ 
In particular, both conditionings of Lemma \ref{reenracinement} make sense. Let $\tr$ be a two-type tree such that $\wh{v}_0\in\dd_0\tr$ and $k_\emptyset(\tr)=1$. Since the trees $\tr$ and $\wh{\tr}^{(\wh{v}_0)}$ represent the same graph, we have
\begin{eqnarray*}
Q_\mu\left(\t^{(\wh{v}_0)}=\tr\right)&=&\prod_{u\in\tr^0\setminus\{\emptyset,\wh{v}_0\}}\mu_0(k_u(\tr))\prod_{u\in\tr^1}\mu_1(k_u(\tr))\\
&=&\prod_{u\in\tr^0\setminus\{\emptyset,\wh{v}_0\}}\mu_0({\rm deg}(u)-1)\prod_{u\in\tr^1}\mu_1({\rm deg}(u)-1)\\
&=&\prod_{u\in\wh{\tr}^{(\wh{v}_0),0}\setminus\{\emptyset,v_0\}}\mu_0({\rm deg}(u)-1)\prod_{u\in\wh{\tr}^{(\wh{v}_0),1}}\mu_1({\rm deg}(u)-1)\\
&=&Q_\mu\left(\t^{(v_0)}=\wh{\tr}^{(\wh{v}_0)}\right)\\
&=&Q_\mu\left(\wh{\t}^{(v_0)}=\tr\right),
\end{eqnarray*}
which implies the desired result.\cq

Before stating a spatial version of Lemma \ref{reenracinement}, we establish a symmetry property of the collection of measures $\nu$. To this end, we let $\wt{\nu}_1^k$ be the image measure of $\nu_1^k$ under the mapping $(x_1,\ldots,x_k)\in\R^k\lmt(x_k,\ldots,x_1)$ and we set $\wt{\nu}=((\nu_0^k,\wt{\nu}_1^k))_{k\geq1}$. 

\begin{lemma}\label{invariance} For every $k\geq1$ and every $j\in\{1,\ldots,k\}$, the measures $\nu_1^k$ and $\wt{\nu}_1^k$ are invariant under the mapping $\phi_j:\R^k\la\R^k$ defined by
$$\phi_j(x_1,\ldots,x_k)=(x_{j+1}-x_j,\ldots,x_k-x_j,-x_j,x_1-x_j,\ldots,x_{j-1}-x_j).$$
\end{lemma}

\proof Recall the definition of the sets $A_k$ and $B_k$. Let $\rho^{k}$ be the uniform distribution on $B_k$. Then $\nu_1^k$ is the image measure of $\rho^k$ under the mapping $\vph_k:B_k\la A_k$ defined by
$$\vph_k(x_1,\ldots,x_{k+1})=(x_1,x_1+x_2,\ldots,x_1+\ldots+x_k).$$
For every $(x_1,\ldots,x_{k+1})\in\R^{k+1}$ we set $p_j(x_1,\ldots,x_{k+1})=(x_{j+1},\ldots,x_{k+1},x_1,\ldots,x_j)$. It is immediate that $\rho^k$ is invariant under the mapping $p_j$. Furthermore $\phi_j\circ\vph_k(x)=\vph_k\circ p_j(x)$ for every $x\in B_k$, which implies that $\nu_1^k$ is invariant under $\phi_j$.  

At last for every $(x_1,\ldots,x_k)\in\R^k$ we set $S(x_1,\ldots,x_k)=(x_k,\ldots,x_1)$. Then $\phi_j\circ S=S\circ\phi_{k-j+1}$, which implies that $\wt{\nu}_1^k$ is invariant under $\phi_j$.  \cq

If $(\t,U)\in\Omega$ and $v_0\in\t^0$, the re-rooted spatial tree $(\wh{\t}^{(v_0)},\wh{U}^{(v_0)})$ is defined as follows. For every vertex $v\in\wh{\t}^{(v_0),0}$, we set 
$$\wh{U}^{(v_0)}_v=U_{\ov{v}}-U_{v_0},$$
where $\ov{v}$ is the vertex of the initial tree $\t$ corresponding to $v$, and for every vertex $v\in\wh{\t}^{(v_0),1}$, we set 
$$\wh{U}^{(v_0)}_v=\wh{U}^{(v_0)}_{\check{v}}.$$
Note that, since $v_0\in\t^0$, $v\in\wh{\t}^{(v_0)}$ is of type $0$ if and only if $\ov{v}\in\t$ is of type $0$.

If $(\t,U)\in\Omega$ and $w_0\in\t$, we also consider the spatial tree $(\t^{(w_0)},U^{(w_0)})$ where $U^{(w_0)}$ is the restriction of $U$ to the tree $\t^{(w_0)}$. 

Recall the definition of the probability measure $\Q_{\mu,\nu}$.

\begin{lemma}\label{reenracinementspatial} Let $v_0\in\u^\ast$ be of the form $v_0=1u^2\ldots u^{2p}$ for some $p\in\N$. Assume that $Q_\mu(v_0\in\t)>0$. Then the law of the re-rooted spatial tree $(\wh{\t}^{(v_0)},\wh{U}^{(v_0)})$ under $\Q_{\mu,\wt{\nu}}(\cdot\mid v_0\in\t)$ coincides with the law of the spatial tree $(\t^{(\wh{v}_0)},U^{(\wh{v}_0)})$ under $\Q_{\mu,\nu}(\cdot\mid\wh{v}_0\in\t)$.
\end{lemma}
Lemma \ref{reenracinementspatial} is a consequence of Lemma \ref{reenracinement} and Lemma \ref{invariance}. We leave details to the reader.

If $(\t,U)\in\Omega$, we denote by $\Delta_0=\Delta_0(\t,U)$ the set of all vertices of type $0$ with minimal spatial position~:
$$\Delta_0=\left\{v\in\t^0:U_v=\min\{U_w:w\in\t\}\right\}.$$
We also denote by $v_m$ the first element of $\Delta_0$ in the lexicographical order. The following two lemmas can be proved from Lemma \ref{reenracinementspatial} in the same way as Lemma 3.3 and Lemma 3.4 in \cite{LG}.

\begin{lemma}\label{reenrmin1}
For any nonnegative measurable functional $F$ on $\Omega$,
$$\Q_{\mu,\wt{\nu}}\left(F\left(\wh{\t}^{(v_m)},\wh{U}^{(v_m)}\right)\ind{\{\#\Delta_0=1,v_m\in\dd_0\t\}}\right)=\Q_{\mu,\nu}\left(F(\t,U)(\#\dd_0\t)\ind{\{\un{U}>0\}}\right).$$
\end{lemma}

\begin{lemma}\label{reenrmin2}
For any nonnegative measurable functional $F$ on $\Omega$,
$$\Q_{\mu,\wt{\nu}}\left(\sum_{v_0\in\Delta_0\cap\dd_0\t}F\left(\wh{\t}^{(v_0)},\wh{U}^{(v_0)}\right)\right)=\Q_{\mu,\nu}\left(F(\t,U)(\#\dd_0\t)\ind{\{\un{U}\geq0\}}\right).$$
\end{lemma}

\subsection{Estimates for the probability of staying on the positive side}

In this section we will derive upper and lower bounds for the probability $\P^n_{\mu,\nu,x}(\un{U}>0)$ as $n\to\infty$. We first state a lemma which is a direct consequence of Lemma 17 in \cite{MaMi}.

\begin{lemma}\label{asympt} There exist constants $c_0>0$ and $c_1>0$ such that
\ba
n^{3/2}P_\mu\left(\#\t^1=n\right)&\build{\la}_{n\to\infty}^{}&c_0,\\
n^{3/2}Q_\mu\left(\#\t^1=n\right)&\build{\la}_{n\to\infty}^{}&c_1.
\ea 
\end{lemma}

We now establish a preliminary estimate concerning the number of leaves of type $0$ in a tree with $n$ vertices of type $1$. 

\begin{lemma}\label{feuilles}
There exists a constant $\beta_0>0$ such that for every $n$ sufficiently large,
$$P_\mu\left(\left|(\#\dd_0\t)-m_1\mu_0(0)n\right|>n^{3/4},\#\t^1=n\right)\leq e^{-n^{\beta_0}}.$$
\end{lemma}

\proof Let $\t$ be a two-type tree. Recall that $\zeta=\#\t-1$. Let 
$$v(0)=\emptyset\prec v(1)\prec\ldots\prec v(\zeta)$$ 
be the vertices of $\t$ listed in lexicographical order. For every $n\in\{0,1,\ldots,\zeta\}$ we define $R_n=(R_n(k))_{k\geq1}$ as follows. For every $k\in\{1,\ldots,|v(n)|\}$, $R_n(k)$ is the number of younger brothers of the ancestor of $v(n)$ at generation $k$. Here younger brothers are those brothers which have not yet been visited at time $n$ in search-depth sequence. For every $k>|v(n)|$, we set $R_n(k)=0$. Standard arguments (see e.g.~\cite{LGLJ} for similar results) show that $(R_n,|v(n)|)_{0\leq n\leq\zeta}$ has the same distribution as $(R'_n,h'_n)_{0\leq n\leq T'-1}$, where $(R'_n,h'_n)_{n\geq0}$ is a Markov chain whose transition kernel is given by~:
\begin{enumerate}
\item[$\bullet$] $S\Big(((r_1,\ldots,r_h,0,\ldots),h),((r_1,\ldots,r_h,k-1,0,\ldots),h+1)\Big)=\mu_i(k)$ for $k\geq1$, $h\geq0$ and $r_1,\ldots,r_h\geq0$,
\item[$\bullet$] $S\Big(((r_1,\ldots,r_h,0,\ldots),h),((r_1,\ldots,r_l-1,0,\ldots),l)\Big)=\mu_i(0),$ where $l=\inf\{m\geq1:r_m>0\}$, for $h\geq1$ and $r_1,\ldots,r_h\geq0$ such that $\{m\geq1:r_m>0\}\neq\emptyset$,
\item[$\bullet$] $S((0,h),(0,0))=\mu_i(0)$ for every $h\geq0$,
\end{enumerate}
where $i=0$ if $h$ is even, and $i=1$ if $h$ is odd, and finally
$$T'=\inf\left\{n\geq1:(R'_n,h'_n)=(0,0)\right\}.$$
Write $\PP'$ for the probability measure under which $(R'_n,h'_n)_{n\geq0}$ is defined. We define a sequence of stopping times $(\tau'_j)_{j\geq0}$ by $\tau'_0=\inf\{n\geq0:h'_n\;{\rm is}\;{\rm odd}\}$ and $\tau'_{j+1}=\inf\{n>\tau'_j:h'_n\;{\rm is}\;{\rm odd}\}$ for every $j\geq0$. At last we set for every $j\geq0$,
$$X'_j=\ind{}\left\{h'_{\tau'_j+1}=h'_{\tau'_j}+1\right\}\left(1+R'_{\tau'_j+1}\left(h'_{\tau'_j+1}\right)\right).$$
Since $\#\t^0=1+\sum_{u\in\t^1}k_u(\t)$, we have
\begin{eqnarray*}
P_\mu\left(|\#\t^0-m_1n|>n^{3/4},\#\t^1=n\right)&=&\PP'\left(\Big|\sum_{j=0}^{n-1}X'_j-m_1n+1\Big|>n^{3/4},\tau'_{n-1}<T'<\tau'_n\right)\\
&\leq&\PP'\left(\Big|\sum_{j=0}^{n-1}X'_j-m_1n\Big|>n^{3/4}-1\right).
\end{eqnarray*}
Thanks to the strong Markov property, the random variables $X'_j$ are independent and distributed according to $\mu_1$. A standard moderate deviations inequality ensures the existence of a positive constant $\beta_1>0$ such that for every $n$ sufficiently large, 
\be\label{beta_1}
P_\mu\left(\left|\#\t^0-m_1n\right|>n^{3/4},\#\t^1=n\right)\leq e^{-n^{\beta_1}}.
\ee
In the same way as previously, we define another sequence of stopping times $(\theta'_j)_{j\geq0}$ by $\theta'_0=0$ and $\theta'_{j+1}=\inf\{n>\theta'_j:h'_n\;{\rm is}\;{\rm even}\}$ for every $j\geq0$ and we set for every $j\geq0$,
$$Y'_j=\ind{}\left\{h'_{\theta'_j+1}\leq h'_{\theta'_j}\right\}.$$
Using the sequences $(\theta'_j)_{j\geq0}$ and $(Y'_j)_{j\geq0}$, an argument similar to the proof of (\ref{beta_1}) shows that there exists a positive constant $\beta_2>0$ such that for every $n$ sufficiently large, 
\be\label{beta_2}
P_\mu\left(\left|\#(\dd_0\t)-\mu_0(0)n\right|>n^{5/8},\#\t^0=n\right)\leq e^{-n^{\beta_2}}.
\ee
From (\ref{beta_1}), we get for $n$ sufficiently large, 
\begin{eqnarray*}
&&\hspace{-1cm}P_\mu\left(\left|(\#\dd_0\t)-m_1\mu_0(0)n\right|>n^{3/4},\#\t^1=n\right)\\
&\leq&e^{-n^{\beta_1}}+P_\mu\left(\left|(\#\dd_0\t)-m_1\mu_0(0)n\right|>n^{3/4},|\#\t^0-m_1n|\leq n^{3/4}\right).\\
\end{eqnarray*}
However, for $n$ sufficiently large,
\begin{eqnarray*}
&&\hspace{-1cm}P_\mu\left(\left|(\#\dd_0\t)-m_1\mu_0(0)n\right|>n^{3/4},|\#\t^0-m_1n|\leq n^{3/4}\right)\\
&=&\sum_{k=\lceil-n^{3/4}+m_1n\rceil}^{\lfloor n^{3/4}+m_1n\rfloor}P_\mu\left(\left|(\#\dd_0\t)-m_1\mu_0(0)n\right|>n^{3/4},\#\t^0=k\right)\\
&\leq&\sum_{k=\lceil-n^{3/4}+m_1n\rceil}^{\lfloor n^{3/4}+m_1n\rfloor}P_\mu\left(\left|(\#\dd_0\t)-\mu_0(0)k\right|>(1-\mu_0(0))n^{3/4},\#\t^0=k\right)\\
&\leq&\sum_{k=\lceil-n^{3/4}+m_1n\rceil}^{\lfloor n^{3/4}+m_1n\rfloor}P_\mu\left(\left|(\#\dd_0\t)-\mu_0(0)k\right|>k^{5/8},\#\t^0=k\right).
\end{eqnarray*}
At last, we use (\ref{beta_2}) to obtain for $n$ sufficiently large,
$$P_\mu\left(\left|(\#\dd_0\t)-m_1\mu_0(0)n\right|>n^{3/4},|\#\t^0-m_1n|\leq n^{3/4}\right)\leq(2n^{3/4}+1)e^{-Cn^{\beta_2}},$$
where $C$ is a positive constant. The desired result follows by combining this last estimate with (\ref{beta_1}).\cq

We will now state a lemma which plays a crucial role in the proof of the main result of this section. To this end, recall the definition of $v_m$ and set for every $n\geq1$,
\ba
Q^n_{\mu}&=&Q_{\mu}\left(\cdot\mid\#\t^1=n\right),\\
\Q^n_{\mu,\nu}&=&\Q_{\mu,\nu}\left(\cdot\mid\#\t^1=n\right).
\ea

\begin{lemma}\label{argmin}
There exists a constant $c>0$ such that for every $n$ sufficiently large,
$$\Q_{\mu,\nu}^n(v_m\in\dd_0\t)\geq c.$$
\end{lemma}

\proof The proof of this lemma is similar to the proof of Lemma 4.3 in \cite{LG}. Nevertheless, we give a few details to explain how this proof can be adapted to our context. 

Choose $p\geq1$ such that $\mu_1(p)>0$. Under $\Q_{\mu,\nu}(\cdot\mid k_1(\t)=p,\;k_{11}(\t)=\ldots=k_{1p}(\t)=2\,)$, we can define $2p$ spatial trees $\{(\t^{ij},U^{ij}),i=1,\ldots,p,j=1,2\}$ as follows. For $i\in\{1,\ldots,p\}$ and $j=1,2$, we set 
$$\t^{ij}=\{\emptyset\}\cup\{1v:1ijv\in\t\},$$
$U^{ij}_\emptyset=0$ and $U^{ij}_{1v}=U_{1ijv}-U_{1i}$ if $1ijv\in\t$. Then under the probability measure $\Q_{\mu,\nu}(\cdot\mid k_1(\t)=p,\;k_{11}(\t)=\ldots=k_{1p}(\t)=2\,)$, the trees $\{(\t^{ij},U^{ij}),i=1,\ldots,p,j=1,2\}$ are independent and distributed according to $\Q_{\mu,\nu}$. Furthermore, we notice that under the measure $\Q_{\mu,\nu}(\cdot\mid k_1(\t)=p,\;k_{11}(\t)=\ldots=k_{1p}(\t)=2\,)$, we have with an obvious notation
\begin{eqnarray*}
&&\hspace{-0.6cm}\Big(\left\{\#\t^{11,1}+\#\t^{12,1}=n-2p+1\right\}\cap\{U_{11}<0\}\cap\left\{v^{11}_m\in\dd_0\t^{11}\right\}\cap\bigcap_{2\leq i\leq p}\left\{U_{1i}\geq0\right\}\\
&&\hspace{-0.4cm}\cap\left\{\un{U}^{12}\geq0\right\}\cap\bigcap_{2\leq i\leq p,j=1,2}\left\{\t^{ij}=\{\emptyset,1,11,\ldots,1p\},\,\un{U}^{ij}\geq0\right\}\Big)\subset\left\{\#\t^1=n,v_m\in\dd_0\t\right\}.
\end{eqnarray*}
So we have for $n\geq1+2p$, 
\begin{eqnarray}
&&\hspace{-1cm}\Q_{\mu,\nu}\left(\#\t^1=n,v_m\in\dd_0\t\right)\nonumber\\
&\geq&\;C(\mu,\nu,p)\sum_{j=1}^{n-2p}\Q_{\mu,\nu}\left(\#\t^1=j,v_m\in\dd_0\t\right)\,\Q_{\mu,\nu}\left(\#\t^1=n-j-2p+1,\,\un{U}\geq0\right)\label{eqargmin},
\end{eqnarray}
where 
$$C(\mu,\nu,p)=\mu_1(p)^{2p-1}\mu_0(2)^p\mu_0(0)^{2p(p-1)}\nu_1^p\left((-\infty,0)\times[0,+\infty)^{p-1}\right)\nu_1^p\left([0,+\infty)^p\right)^{2(p-1)}.$$ 
From (\ref{eqargmin}), we are now able to get the result by following the lines of the proof of Lemma 4.3 in \cite{LG}.\cq

We can now state the main result of this section.

\begin{proposition}\label{estpositif} Let $K>0$. There exist constants $\gamma_1>0$, $\gamma_2>0$, $\wt{\gamma}_1>0$ and $\wt{\gamma}_2>0$ such that for every $n$ sufficiently large and for every $x\in[0,K]$,
$$\frac{\wt{\gamma}_1}{n}\leq\Q_{\mu,\nu}^n(\,\un{U}>0)\leq\frac{\wt{\gamma}_2}{n},$$
$$\frac{\gamma_1}{n}\leq\P_{\mu,\nu,x}^n(\,\un{U}>0)\leq\frac{\gamma_2}{n}.$$
\end{proposition}

\proof The proof of Proposition \ref{estpositif} is similar to the proof of Proposition 4.2 in \cite{LG}. The major difference comes from the fact that we cannot easily get an upper bound for $\#(\dd_0\t)$ on the event $\{\#\t^1=n\}$. In what follows, we will explain how to circumvent this difficulty.

We first use Lemma \ref{reenrmin1} with $F(\t,U)=\ind{\{\t^1=n\}}$. Since $\#\wh{\t}^{(v_0),1}=\#\t^1$ if $v_0\in\dd_0\t$, we get
\be\label{eqestpos0}
\Q_{\mu,\nu}\left(\#(\dd_0\t)\ind{\{\#\t^1=n,\,\un{U}>0\}}\right)\leq Q_\mu\left(\#\t^1=n\right).
\ee
On the other hand, we have
\begin{eqnarray}
&&\hspace{-1cm}\Q_{\mu,\nu}\left(\#(\dd_0\t)\ind{\{\#\t^1=n,\,\un{U}>0\}}\right)\nonumber\\
&\geq&\;m_1\mu_0(0)n\Q_{\mu,\nu}\left(\#\t^1=n,\,\un{U}>0\right)-\Q_{\mu,\nu}\left(|\#(\dd_0\t)-m_1\mu_0(0)n|\ind{\{\#\t^1=n,\,\un{U}>0\}}\right)\label{eqestpos1}.
\end{eqnarray}
Now thanks to Lemma \ref{feuilles}, we have for $n$ sufficiently large,
\begin{eqnarray}
&&\hspace{-1cm}\Q_{\mu,\nu}\left(|\#(\dd_0\t)-m_1\mu_0(0)n|\ind{\{\#\t^1=n,\,\un{U}>0\}}\right)\nonumber\\
&\leq&\;n^{3/4}\,\Q_{\mu,\nu}\left(\#\t^1=n,\,\un{U}>0\right)+\Q_{\mu,\nu}\left(\#(\dd_0\t)\ind{\{\#\t^1=n,\,\un{U}>0\}}\right)\nonumber\\
&&+\;m_1\mu_0(0)n\Q_{\mu,\nu}(|\#(\dd_0\t)-m_1\mu_0(0)n|>n^{3/4},\#\t^1=n)\nonumber\\
&\leq&\;n^{3/4}\,\Q_{\mu,\nu}\left(\#\t^1=n,\,\un{U}>0\right)+\Q_{\mu,\nu}\left(\#(\dd_0\t)\ind{\{\#\t^1=n,\,\un{U}>0\}}\right)+\;m_1\mu_0(0)ne^{-n^{\beta_0}}\label{eqestpos2}.
\end{eqnarray}
From (\ref{eqestpos0}), (\ref{eqestpos1}) and (\ref{eqestpos2}) we get for $n$ sufficiently large
$$(m_1\mu_0(0)n-n^{3/4})\Q_{\mu,\nu}\left(\#\t^1=n,\,\un{U}>0\right)\leq2\,Q_\mu\left(\#\t^1=n\right)+m_1\mu_0(0)ne^{-n^{\beta_0}}.$$
Using Lemma \ref{asympt} it follows that
$$\limsup_{n\to\infty}n\Q_{\mu,\nu}^n\left(\un{U}>0\right)\leq\frac{2}{m_1\mu_0(0)},$$
which ensures the existence of $\wt{\gamma}_2$.

Let us now use Lemma \ref{reenrmin2} with 
$$F(\t,U)=\ind{\{\t^{1}=n,\#(\dd_0\t)\leq m_1\mu_0(0)n+n^{3/4}\}}.$$ 
Since $\#(\dd_0\t)=\#(\dd_0\wh{\t}^{(v_0)})$ if $v_0\in\dd_0\t$, we have for $n$ sufficiently large,
\begin{eqnarray}
&&\Q_{\mu,\nu}\left(\#(\dd_0\t)\ind{\{\#\t^{1}=n,\,\#(\dd_0\t)\leq m_1\mu_0(0)+n^{3/4},\,\un{U}\geq0\}}\right)\nonumber\\
&&\hspace{0.6cm}=\;\Q_{\mu,\wt{\nu}}\left(\#(\Delta_0\cap\dd_0\t)\ind{\{\#\t^{1}=n,\,\#(\dd_0\t)\leq m_1\mu_0(0)+n^{3/4}\}}\right)\nonumber\\
&&\hspace{0.6cm}=\;\Q_{\mu,\nu}\left(\#(\Delta_0\cap\dd_0\t)\ind{\{\#\t^{1}=n,\,\#(\dd_0\t)\leq m_1\mu_0(0)+n^{3/4}\}}\right)\nonumber\\
&&\hspace{0.6cm}\geq\;\Q_{\mu,\nu}\left(\#(\Delta_0\cap\dd_0\t)\geq1,\,\#\t^{1}=n,\,\#(\dd_0\t)\leq m_1\mu_0(0)+n^{3/4}\right)\nonumber\\
&&\hspace{0.6cm}\geq\;\Q_{\mu,\nu}\left(v_m\in\dd_0\t,\,\#\t^{1}=n,\,\#(\dd_0\t)\leq m_1\mu_0(0)+n^{3/4}\right)\nonumber\\
&&\hspace{0.6cm}\geq\;\Q_{\mu,\nu}\left(v_m\in\dd_0\t,\,\#\t^{1}=n\right)-\Q_{\mu,\nu}\left(\#(\dd_0\t)>m_1\mu_0(0)+n^{3/4},\#\t^1=n\right)\nonumber\\
&&\hspace{0.6cm}\geq\;c\,\Q_{\mu,\nu}\left(\#\t^{1}=n\right)-e^{-n^{\beta_0}}\label{eqestpos3},
\end{eqnarray}
where the last inequality comes from Lemma \ref{feuilles} and Lemma \ref{argmin}. On the other hand,
\begin{eqnarray}
&&\Q_{\mu,\nu}\left(\#(\dd_0\t)\ind{\{\#\t^{1}=n,\,\#(\dd_0\t)\leq m_1\mu_0(0)+n^{3/4},\,\un{U}\geq0\}}\right)\nonumber\\
&&\hspace{0.6cm}\leq\;\left(m_1\mu_0(0)n+n^{3/4}\right)\Q_{\mu,\nu}\left(\#\t^{1}=n,\,\un{U}\geq0\right).\label{eqestpos4}
\end{eqnarray}
Then (\ref{eqestpos3}), (\ref{eqestpos4}) and Lemma \ref{asympt} imply that for $n$ sufficiently large,
\be\label{eqestpos5}
\liminf_{n\to\infty}n\Q_{\mu,\nu}^n(\un{U}\geq0)\geq\frac{c}{m_1\mu_0(0)}.
\ee
Recall that $p\geq1$ is such that $\mu_1(p)>0$. Also recall the definition of the spatial tree $(\t^{[w_0]},U^{[w_0]})$. From the proof of Corallary \ref{corconvQ}, we have
\begin{eqnarray}
\P_{\mu,\nu,0}\left(\,\un{U}>0,\#\t^{1}=n\right)&\geq&\P_{\mu,\nu,0}\Big(k_\emptyset(\t)=1,k_1(\t)=p,U_{11}>0,\ldots,U_{1p}>0\nonumber,\\
&&\hspace{1.5cm}\#\t^{[11],1}=n-1,\,\un{U}^{[11]}\geq0,\t^{[12]}=\ldots=\t^{[1p]}=\{\emptyset\}\Big)\nonumber\\
&\geq&C_2(\mu,\nu,p)\P_{\mu,\nu,0}\left(\,\un{U}\geq0,\#\t^{1}=n-1\right)\nonumber\\
&\geq&\mu_0(1)C_2(\mu,\nu,p)\Q_{\mu,\nu}\left(\,\un{U}\geq0,\#\t^{1}=n-1\right)\label{eqestpos6},
\end{eqnarray}
where we have set 
$$C_2(\mu,\nu,p)=\mu_0(1)\mu_1(p)\mu_0(0)^{p-1}\nu_1^p\left((0,+\infty)^p\right).$$ 
We then deduce the existence of $\gamma_1$ from (\ref{eqestpos5}), (\ref{eqestpos6}) and Lemma \ref{asympt}. A similar argument gives the existence of $\wt{\gamma}_1$. 

At last, we define $m\in\N$ by the condition $(m-1)p<K\leq mp$. For every $l\in\N$, we define $1^l\in\u$ by $1^l=11\ldots1$, $|1^l|=l$. Notice that $\nu_1^p(\{(p,p-1,\ldots,1)\})=N(p+1)^{-1}$. By arguing on the event 
$$\left\{k_\emptyset(\t)=k_{11}(\t)=\ldots=k_{1^{2m-2}}=1,k_1(\t)=\ldots=k_{1^{2m-1}}(\t)=p\right\},$$ 
we see that for every $n\geq m$,
\be\label{eqestpos7}
\Q_{\mu,\nu}\left(\#\t^{1}=n,\,\un{U}>0\right)\geq C_3(\mu,\nu,p,m)\P_{\mu,\nu,K}\left(\#\t^{1}=n-m,\,\un{U}>0\right),
\ee
with 
$$C_3(\mu,\nu,p,m)=\mu_0(1)^{m-1}\mu_1(p)^m\mu_0(0)^{m(p-1)}N(p+1)^{-m}.$$ Thanks to Lemma \ref{asympt}, (\ref{eqestpos7}) yields for every $n$ sufficiently large,
$$\P_{\mu,\nu,K}^n(\un{U}>0)\leq\frac{2c_1}{c_0}\,\frac{\wt{\gamma}_2}{C_3(\mu,\nu,p,m)}\,\frac{1}{n},$$
which gives the existence of $\gamma_2$.\cq

\subsection{Asymptotic properties of conditioned trees}

We first introduce a specific notation for rescaled contour and spatial contour processes. For every $n\geq1$ and every $t\in[0,1]$, we set
\ba
C^{(n)}(t)&=&\frac{\sqrt{\rho_\q(Z_\q-1)}}{4}\;\frac{C(2(\#\t-1)t)}{n^{1/2}},\\
V^{(n)}(t)&=&\left(\frac{9(Z_\q-1)}{4\rho_\q}\right)^{1/4}\;\frac{V(2(\#\t-1)t)}{n^{1/4}}.
\ea
In this section, we will get some information about asymptotic properties of the pair $(C^{(n)},V^{(n)})$ under $\ov{\P}^{\,n}_{\mu,\nu,x}$. We will consider the conditioned measure
$$\ov{\Q}^{\,n}_{\mu,\nu}=\Q^n_{\mu,\nu}(\cdot\mid\un{U}>0).$$
Berfore stating the main result of this section, we will establish three lemmas. The first one is the analogue of Lemma 6.2 in \cite{LG} for two-type spatial trees and can be proved in a very similar way.

\begin{lemma}\label{majQ}
There exists a constant $\ov{c}>0$ such that, for every measurable function $F$ on $\Omega$ with $0\leq F\leq 1$,
$$\ov{\Q}_{\mu,\nu}^{\,n}(F(\t,U))\leq\ov{c}\,\Q_{\mu,\wt{\nu}}^{n}\left(F\left(\wh{\t}^{(v_{m})},\wh{U}^{(v_{m})}\right)\right)+O\left(n^{5/2}e^{-n^{\beta_0}}\right),$$
where the constant $\beta_0$ is defined in Lemma \ref{feuilles} and the estimate $O\left(n^{5/2}e^{-n^{\beta_0}}\right)$ for the remainder holds uniformly in $F$.
\end{lemma}

Recall the notation $\check{v}$ for the ``father'' of the vertex $v\in\t\setminus\{\emptyset\}$. 

\begin{lemma}\label{dperefils}
For every $\vep>0$, 
\be\label{dperefils1}
\Q_{\mu,\nu}^{n}\left(\sup_{v\in\t\setminus\{\emptyset\}}\frac{|U_v-U_{\check{v}}|}{n^{1/4}}>\vep\right)\build{\la}_{n\to\infty}^{}0.
\ee
Likewise, for every $\vep>0$ and $x\geq0$, 
\be\label{dperefils2}
\ov{\P}_{\mu,\nu,x}^{\,n}\left(\sup_{v\in\t\setminus\{\emptyset\}}\frac{|U_v-U_{\check{v}}|}{n^{1/4}}>\vep\right)\build{\la}_{n\to\infty}^{}0.
\ee
\end{lemma}

\proof Let $\vep>0$. First notice that the probability measure $\nu_1^k$ is supported on the set $\{-k,-k+1,\ldots,k\}^k$. Then we have $\Q_{\mu,\nu}^n$ a.s.~or $\ov{\P}^{\,n}_{\mu,\nu,x}$ a.s.,
$$\sup_{v\in\t\setminus\{\emptyset\}}|U_v-U_{\check{v}}|=\sup_{v\in\t^0\setminus\{\emptyset\}}|U_v-U_{\check{v}}|\leq\sup_{v\in\t^1}k_v(\t).$$
Now, from Lemma 16 in \cite{MaMi}, there exists a constant $\al_0>0$ such that for all $n$ sufficiently large,
\ba
Q_\mu\left(\sup_{v\in\t^1}k_v(\t)>\vep n^{1/4}\right)&\leq&e^{-n^{\al_0}},\\
P_\mu\left(\sup_{v\in\t^1}k_v(\t)>\vep n^{1/4}\right)&\leq&e^{-n^{\al_0}}.
\ea
Our assertions (\ref{dperefils1}) and (\ref{dperefils2}) easily follow using also Lemma \ref{asympt} and Proposition \ref{estpositif}.\cq
 
Recall the definition of the re-rooted tree $(\wh{\t}^{(v_0)},\wh{U}^{(v_0)})$. Its contour and spatial contour functions $(\wh{C}^{(v_0)},\wh{V}^{(v_0)})$ are defined on the line interval $[0,2(\#\wh{\t}^{(v_0)}-1)]$. We extend these functions to the line interval $[0,2(\#\t-1)]$ by setting $\wh{C}^{(v_0)}(t)=0$ and $\wh{V}^{(v_0)}(t)=0$ for every $t\in[2(\#\wh{\t}^{(v_0)}-1),2(\#\t-1)]$. Also recall the definition of $v_m$.

At last, recall that $(\ov{\e}^0,\ov{\r}^0)$ denotes the conditioned Brownian snake.

\begin{lemma}\label{convreenr}
The law under $\Q^n_{\mu,\nu}$ of
$$\left(\left(\frac{\sqrt{\rho_\q(Z_\q-1)}}{4}\;\frac{\wh{C}^{(v_m)}(2(\#\t-1)t)}{n^{1/2}}\right)_{0\leq t\leq1},\left(\left(\frac{9(Z_\q-1)}{4\rho_\q}\right)^{1/4}\frac{\wh{V}^{(v_m)}(2(\#\t-1)t)}{n^{1/4}}\right)_{0\leq t\leq1}\right)$$
converges to the law of $(\ov{\e}^0,\ov{\r}^0)$. The convergence holds in the sense of weak convergence of probability measures on the space $C([0,1],\R)^2$. 
\end{lemma}

\proof From Corollary \ref{corconvQ} and the Skorokhod representation theorem, we can construct on a suitable probability space a sequence  a sequence $(\t_n,U_n)$ and a Brownian snake $(\eb,\rb^0)$, such that each pair $(\t_n,U_n)$ is distributed according to $\Q^n_{\mu,\nu}$, and such that if we write $(C_n,V_n)$ for the contour and spatial contour functions of $(\t_n,U_n)$ and $\zeta_n=\#\t_n-1$, we have
\be\label{convskho}
\left(\left(\frac{\sqrt{\rho_\q(Z_\q-1)}}{4}\;\frac{C_n(2\zeta_nt)}{n^{1/2}}\right)_{0\leq t\leq1},\left(\left(\frac{9(Z_\q-1)}{4\rho_\q}\right)^{1/4}\frac{V_n(2\zeta_nt)}{n^{1/4}}\right)_{0\leq t\leq1}\right)\build{\la}_{n\to\infty}^{}(\eb,\rb^0),
\ee
uniformly on $[0,1]$, a.s. 

Then if $(\t,U)\in\Omega$ and $v_0\in\t$, we introduce a new spatial tree $(\wh{\t}^{(v_0)},\wt{U}^{(v_0)})$ by setting for every $w\in\wh{\t}^{(v_0)}$ 
$$\wt{U}^{(v_0)}_w=U_{\ov{w}}-U_{v_0},$$ 
where $\ov{w}$ is the vertex corresponding to $w$ in the initial tree (in contrast with the definition of $\wh{U}^{(v_0)}_w$, $\wt{U}^{(v_0)}_v$ does not necesseraly coincide with $\wt{U}^{(v_0)}_{\check{v}}$ when $v$ is of type $1$). We denote by $\wt{V}^{(v_0)}$ the spatial contour function of $(\wh{\t}^{(v_0)},\wt{U}^{(v_0)})$, and we set $\wt{V}^{(v_0)}(t)=0$ for $t\in[2(\#\wh{\t}^{(v_0)}-1),2(\#\t-1)]$. Note that, if $w\in\wh{\t}^{(v_0)}$ is either a vertex of type $0$ or a vertex of type $1$ which does not belong to the ancestral line of $\wh{v}_0$, then 
$$\wh{U}^{(v_0)}_w=\wt{U}^{(v_0)}_w,$$ 
whereas if $w\in\wh{\t}^{(v_0),1}$ belongs to the ancestral line of $\wh{v}_0$, then 
$$\wh{U}^{(v_0)}_w=\wt{U}^{(v_0)}_{\check{w}}.$$ 
Then we have
\be\label{diffwhwtU}
\sup_{w\in\wh{\t}^{(v_0)}}\;\left|\wh{U}_{w}-\wt{U}_{w}\right|\leq\sup_{w\in\t\setminus\{\emptyset\}}\;\left|U_{w}-U_{\check{w}}\right|.
\ee

Write $v_m^n$ for the first vertex realizing the minimal spatial position in $\t_n$. In the same way as in the derivation of (18) in the proof of Proposition 6.1 in \cite{LG}, it follows from (\ref{convskho}) that
$$\left(\left(\frac{\sqrt{\rho_\q(Z_\q-1)}}{4}\;\frac{\wh{C}^{(v_m^n)}_n(2\zeta_nt)}{n^{1/2}}\right)_{0\leq t\leq1},\left(\left(\frac{9(Z_\q-1)}{4\rho_\q}\right)^{1/4}\frac{\wt{V}^{(v_m^n)}_n(2\zeta_nt)}{n^{1/4}}\right)_{0\leq t\leq1}\right)\build{\la}_{n\to\infty}^{}(\ov{\eb}^0,\ov{\rb}^0),$$
uniformly on $[0,1]$, a.s., where the conditioned pair $(\ov{\eb}^0,\ov{\rb}^0)$ is constructed from the unconditioned one $(\eb,\rb^0)$ as explained in section \ref{serpent}. Let $\vep>0$. We deduce from (\ref{diffwhwtU}) that
$$\PP'\left(\sup_{t\in[0,1]}\left|\frac{\wh{V}^{(v_m^n)}_n(2\zeta_nt)}{n^{1/4}}-\frac{\wt{V}^{(v_m^n)}_n(2\zeta_nt)}{n^{1/4}}\right|>\vep\right)\leq\Q_{\mu,\nu}^n\left(\sup_{w\in\t\setminus\{\emptyset\}}\;\frac{\left|U_{w}-U_{\check{w}}\right|}{n^{1/4}}>\vep\right),$$
where we have written $\PP'$ for the probability measure under which the sequence $(\t_n,U_n)_{n\geq1}$ and the Brownian snake $(\eb,\rb^0)$ are defined. From (\ref{dperefils1}) we get
$$\PP'\left(\sup_{t\in[0,1]}\left|\frac{\wh{V}^{(v_m^n)}_n(2\zeta_nt)}{n^{1/4}}-\frac{\wt{V}^{(v_m^n)}_n(2\zeta_nt)}{n^{1/4}}\right|>\vep\right)\build{\la}_{n\to\infty}^{}0,$$
and the desired result follows.\cq

The following proposition can be proved using Lemma \ref{majQ} and Lemma \ref{convreenr} in the same way as Proposition 6.1 in \cite{LG}.
\begin{proposition}\label{asymptQP} 
For every $b>0$ and $\vep\in(0,1/10)$, we can find $\al,\dl\in(0,\vep)$ such that for all $n$ sufficiently large,
$$\ov{\Q}^{\,n}_{\mu,\nu}\left(\inf_{t\in[\dl/2,1-\dl/2]}V^{(n)}(t)\geq2\al,\sup_{t\in[0,4\dl]\cup[1-4\dl]}\left(C^{(n)}(t)+V^{(n)}(t)\right)\leq\frac{\vep}{2}\right)\geq1-b.$$
Consequently, if $K>0$, we have also for all $n$ sufficiently large, for every $x\in[0,K]$,
$$\ov{\P}^{\,n}_{\mu,\nu,x}\left(\inf_{t\in[\dl,1-\dl]}V^{(n)}(t)\geq\al,\sup_{t\in[0,3\dl]\cup[1-3\dl]}\left(C^{(n)}(t)+V^{(n)}(t)\right)\leq\vep\right)\geq1-\gamma_{3}b,$$
where the constant $\gamma_{3}$ only depends on $\mu,\nu,K$.
\end{proposition}

\subsection{Proof of Theorem \ref{thconvcond}}

The proof below is similar to Section 7 in \cite{LG}. We provide details because the fact that we deal with two-type trees creates nontrivial additional difficulties. 

On a suitable probability space $(\OO,\PP)$ we can define a collection of processes $(\e,\r^z)_{z\geq0}$ such that $(\e,\r^z)$ is a Brownian snake with initial point $z$ for every $z\geq0$. Recall from section \ref{serpent} the definition of $(\ov{\e}^z,\ov{\r}^z)$ and the construction of the conditioned Brownian snake $(\ov{\e}^0,\ov{\r}^0)$.

Recall that $C([0,1],\R)^2$ is equipped with the norm $\|(f,g)\|=\|f\|_u\vee\|g\|_u$. For every $f\in C([0,1],\R)$ and $r>0$, we set
$$\omega_f(r)=\sup_{s,t\in[0,1],|t-s|\leq r}|f(s)-f(t)|.$$

Let $x\geq0$ be fixed throughout this section and let $F$ be a bounded Lipschitz function. We have to prove that
$$\ov{\E}_{\mu,\nu,x}^{\,n}\left(F\left(C^{(n)},V^{(n)}\right)\right)\build{\la}_{n\to\infty}^{}\EE\left(F(\ov{\e}^0,\ov{\r}^0)\right).$$
We may and will assume that $0\leq F\leq1$ and that the Lipschitz constant of $F$ is less than $1$. 

The first lemma we have to prove gives a spatial Markov property for our spatial trees. We use  the notation of section 5 in \cite{LG}. Let recall briefly this notation. We fix $a>0$. If $(\t,U)$ is a mobile and $v\in\t$, we say that $v$ is an exit vertex from $(-\infty,a)$ if $U_v\geq a$ and $U_{v'}<a$ for every ancestor $v'$ of $v$ distinct from $v$. Notice that, since $U_v=U_{\check{v}}$ for every $v\in\t^1$, an exit vertex is necessarily of type $0$. We denote by $v_1,\ldots,v_M$ the exit vertices from $(-\infty,a)$ listed in lexicographical order. For $v\in\t$, recall that $\t^{[v]}=\{w\in\u:vw\in\t\}$. For every $w\in\t^{[v]}$ we set
$$\ov{U}^{[v]}_w=U_{vw}=U^{[v]}_w+U_v.$$  
At last, we denote by $\t^a$ the subtree of $\t$ consisting of those vertices which are not strict descendants of $v_1,\ldots,v_M$. Note in particular that $v_1,\ldots,v_M\in\t^a$. We also write $U^a$ for the restriction of $U$ to $\t^a$. The tree $(\t^a,U^a)$ corresponds to the tree $(\t,U)$ which has been truncated at the first exit time from $(-\infty,a)$. The following lemma is an easy application of classical properties of Galton-Watson trees. We leave details of the proof to the reader.

\begin{lemma}\label{Markovsp} Let $x\in[0,a)$ and $p\in\{1,\ldots,n\}$. Let $n_1,\ldots,n_p$ be positive integers such that $n_1+\ldots+n_p\leq n$. Assume that
$$\ov{\P}_{\mu,\nu,x}^{\,n}\left(M=p,\;\#\t^{[v_1],1}=n_1,\ldots,\;\#\t^{[v_p],1}=n_p\right)>0.$$
Then, under the probability measure $\ov{\P}_{\mu,\nu,x}^{\,n}(\cdot\mid M=p,\,\#\t^{[v_1],1}=n_1,\ldots,\,\#\t^{[v_p],1}=n_p)$, and conditionally on $(\t^a,U^a)$, the spatial trees 
$$\left(\t^{[v_1]},\ov{U}^{[v_1]}\right),\ldots,\left(\t^{[v_p]},\ov{U}^{[v_p]}\right)$$ 
are independent and distributed respectively according to $\ov{\P}_{\mu,\nu,U_{v_1}}^{\,n},\ldots,\ov{\P}_{\mu,\nu,U_{v_p}}^{\,n}$.
\end{lemma}

The next lemma is analogous to Lemma 7.1 in \cite{LG} and can be proved in the same way using Theorem \ref{thconv}.

\begin{lemma}\label{condpos}
Let $0<c'<c''$. Then
$$\sup_{c'n^{1/4}\leq y\leq c''n^{1/4}}\left|\ov{\E}_{\mu,\nu,y}^{\,n}\left(F\left(C^{(n)},V^{(n)}\right)\right)-\EE\left(F\left(\ov{\e}^{By/n^{1/4}},\ov{\r}^{By/n^{1/4}}\right)\right)\right|\build{\la}_{n\to\infty}^{}0,$$
where $B=(9(Z_\q-1)/(4\rho_\q))^{1/4}$.
\end{lemma}

We can now follow the lines of section 7 in \cite{LG}. Let $b>0$. We will prove that for $n$ sufficiently large,
$$\left|\ov{\E}_{\mu,\nu,x}^{\,n}\left(F\left(C^{(n)},V^{(n)}\right)\right)-\EE\left(F\left(\ov{\e}^0,\ov{\r}^0\right)\right)\right|\leq17b.$$ 
We can choose $\vep\in(0,b\wedge1/10)$ in such a way that
\be\label{contserp}
\left|\EE\left(F\left(\ov{\e}^z,\ov{\r}^z\right)\right)-\EE\left(F\left(\ov{\e}^0,\ov{\r}^0\right)\right)\right|\leq b,
\ee
for every $z\in(0,2\vep)$. By taking $\vep$ smaller if necessary, we may also assume that, 
\begin{eqnarray}
\EE\left(\left(3\vep\sup_{0\leq t\leq1}\ov{\e}^0(t)\right)\wedge1\right)\leq b,&&\EE\left(\omega_{\ov{\e}^0}(6\vep)\wedge1\right)\leq b,\nonumber\\
\label{veppetit}
\EE\left(\left(3\vep\sup_{0\leq t\leq1}\ov{\r}^0(t)\right)\wedge1\right)\leq b,&&\EE\left(\omega_{\ov{\r}^0}(6\vep)\wedge1\right)\leq b.
\end{eqnarray}
For $\al,\dl>0$, we denote by $\Gamma_n=\Gamma_n(\al,\dl)$ the event
$$\Gamma_n=\left\{\inf_{t\in[\dl,1-\dl]}V^{(n)}(t)\geq\al,\sup_{t\in[0,3\dl]\cup[1-3\dl,1]}\left(C^{(n)}(t)+V^{(n)}(t)\right)\leq\vep\right\}.$$
From Proposition \ref{asymptQP}, we may fix $\al,\dl\in(0,\vep)$ such that, for all $n$ sufficiently large,
\be\label{P(Gamma_n)}
\ov{\P}_{\mu,\nu,x}^{\,n}(\Gamma_n)>1-b.
\ee
We also require that $\dl$ satisfies the following bound
\be\label{dlm_1}
4\dl(m_1+1)<3\vep.
\ee
Recall the notation $\zeta=\#\t-1$ and $B=(9(Z_\q-1)/(4\rho_\q))^{1/4}$. On the event $\Gamma_n$, we have for every $t\in[2\zeta\dl,2\zeta(1-\dl)]$,
\be\label{Gamma1}
V(t)\ge\al B^{-1}\,n^{1/4}.
\ee

This incites us to apply Lemma \ref{Markovsp} with $a_n=\ov{\al}n^{1/4}$, where $\ov{\al}=\al B^{-1}$. Once again, we use the notation of \cite{LG}. We write $v_1^n,\ldots,v_{M_n}^n$ for the exit vertices from $(-\infty,\ov{\al}n^{1/4})$ of the spatial tree $(\t,U)$, listed in lexicographical order. Consider the spatial trees 
$$\left(\t^{[v_1^n]},\ov{U}^{[v_1^n]}\right),\ldots,\left(\t^{[v_{M_n}^n]},\ov{U}^{[v_{M_n}^n]}\right).$$ 
The contour functions of these spatial trees can be obtained in the following way. Set
$$k_1^n=\inf\left\{k\geq0:V(k)\geq\ov{\al}n^{1/4}\right\},\;\;\;l_1^n=\inf\left\{k\geq k_1^n:C(k+1)<C(k_1^n)\right\},$$
and by induction on $i$,
$$k_{i+1}^n=\inf\left\{k>l_i^n:V(k)\geq\ov{\al}n^{1/4}\right\},\;\;\;l_{i+1}^n=\inf\left\{k\geq k_{i+1}^n:C(k+1)<C(k_i^n)\right\}.$$
Then $k_i^n\leq l_i^n<\infty$ if and only if $i\leq M_n$. Furthermore, $(C(k_i^n+t)-C(k_i^n),0\leq t\leq l_i^n-k_i^n)$ is the contour function of $\t^{[v_i^n]}$ and $(V(k_i^n+t),0\leq t\leq l_i^n-k_i^n)$ is the spatial contour function of $(\t^{[v_i^n]},\ov{U}^{[v_i^n]})$. Using (\ref{Gamma1}), we see that on the event $\Gamma_n$, all integer points of $[2\zeta\dl,2\zeta(1-\dl)]$ must be contained in a single interval $[k_i^n,l_i^n]$, so that for this particular interval we have 
$$l_i^n-k_i^n\geq2\zeta(1-\dl)-2\zeta\dl-2\geq2\zeta(1-3\dl),$$
if $n$ is sufficiently large, $\ov{\P}^{\,n}_{\mu,\nu,x}$ a.s. Hence if
$$E_n=\left\{\exists i\in\{1,\ldots,M_n\}:l_i^n-k_i^n>2\zeta(1-3\dl)\right\},$$
then, for all $n$ sufficiently large, $\Gamma_n\subset E_n$ so that 
\be\label{P(E_n)}
\ov{\P}^{\,n}_{\mu,\nu,x}(E_n)>1-b.
\ee
As in \cite{LG}, on the event $E_n$, we denote by $i_n$ the unique integer $i\in\{1,\ldots,M_n\}$ such that $l_i^n-k_i^n>2\zeta(1-3\dl)$. We also define $\zeta_n=\#\t^{[v_{i_n}^n]}-1$ and $Y_n=U_{v_{i_n}^n}$. Note that $\zeta_n=(l_i^n-k_i^n)/2$ so that $\zeta_n>\zeta(1-3\dl)$. Furthermore, we set
$$p_n=\#\t^{[v_{i_n}^n],1}.$$
We need to prove a lemma providing an estimate of the probability for $p_n$ to be close to $n$. Note that $p_n\leq n=\#\t^1$, $\ov{\P}_{\mu,\nu,x}^{\,n}$ a.s. Recall that $m_1$ denotes the mean of $\mu_1$.

\begin{lemma}\label{p_n}
For every $n$ sufficiently large,
$$\ov{\P}_{\mu,\nu,x}^{\,n}\left(\Gamma_n\cap\{p_n\geq(1-4\dl(m_1+1))n\}\right)\geq1-2b.$$
\end{lemma}

\proof In the same way as in the proof of the bound (\ref{beta_1}) we can verify that there exists a constant $\beta_3>0$ such that for all $n$ sufficiently large, 
$$P_\mu\left(|\zeta-(m_1+1)n|>n^{3/4},\#\t^1=n\right)\leq e^{-n^{\beta_3}}.$$ 
So Lemma \ref{asympt} and Proposition \ref{estpositif} imply that for all $n$ sufficiently large, 
\be\label{eqp_n}
\ov{\P}_{\mu,\nu,x}^{\,n}\left(|\zeta-(m_1+1)n|>n^{3/4}\right)\leq b.
\ee
Now, on the event $\Gamma_n$, we have
$$n-p_n=\#\left(\left(\t\setminus\t^{[v_{i_n}^n]}\right)\cap\t^1\right)\leq\#\left(\t\setminus\t^{[v_{i_n}^n]}\right)=\zeta-\zeta_n\leq3\dl\zeta,$$
since we saw that $\Gamma_n\subset E_n$ and that $\zeta_n>(1-3\dl)\zeta$ on $E_n$. If $n$ is sufficiently large, we have $3\dl(m_1+1)n+3\dl n^{3/4}\leq4\dl(m_1+1)n$ so we obtain that
$$\left(\Gamma_n\cap\left\{\zeta\leq(m_1+1)n+n^{3/4}\right\}\right)\subset\left(\Gamma_n\cap\{p_n\geq(1-4\dl(m_1+1))n\}\right),$$
for all $n$ sufficiently large. The desired result then follows from (\ref{P(Gamma_n)}) and (\ref{eqp_n}).\cq

Let us now define on the event $E_n$, for every $t\in[0,1]$,
\ba
\wt{C}^{(n)}(t)&=&\frac{\sqrt{\rho_q(Z_q-1)}}{4}\;\frac{C\left(k^n_{i_n}+2\zeta_nt\right)-C\left(k^n_{i_n}\right)}{p_n^{1/2}},\\
\wt{V}^{(n)}(t)&=&\left(\frac{9(Z_q-1)}{4\rho_q}\right)^{1/4}\;\frac{V\left(k^n_{i_n}+2\zeta_nt\right)}{p_n^{1/4}}.
\ea
Note that $\wt{C}^{(n)}$ and $\wt{V}^{(n)}$ are rescaled versions of the contour and the spatial contour functions of $(\t^{[v^n_{i_n}]},\ov{U}^{[v^n_{i_n}]})$. On the event $E_n^c$, we take $\wt{C}^{(n)}(t)=\wt{V}^{(n)}(t)=0$ for every $t\in[0,1]$. Straightforward calculations show that on the event $\Gamma_n$, for every $t\in[0,1]$,
\ba
\left|C^{(n)}(t)-\wt{C}^{(n)}(t)\right|&\leq&\vep+\left(1-\frac{p_n^{1/2}}{n^{1/2}}\right)\sup_{s\in[0,1]}\wt{C}^{(n)}(s)+\omega_{\wt{C}^{(n)}}(6\dl),\\
\left|V^{(n)}(t)-\wt{V}^{(n)}(t)\right|&\leq&\vep+\left(1-\frac{p_n^{1/4}}{n^{1/4}}\right)\sup_{s\in[0,1]}\wt{V}^{(n)}(s)+\omega_{\wt{V}^{(n)}}(6\dl).
\ea
Set 
$$\wt{\Gamma}_n=\Gamma_n\cap\{p_n\geq(1-4\dl(m_1+1))n\}.$$ 
We then get that on the event $\wt{\Gamma}_n$, for every $t\in[0,1]$,
\begin{eqnarray}
\label{majwtGamma1}\left|C^{(n)}(t)-\wt{C}^{(n)}(t)\right|&\leq&\vep+4\dl(m_1+1)\sup_{s\in[0,1]}\wt{C}^{(n)}(s)+\omega_{\wt{C}^{(n)}}(6\dl),\\
\label{majwtGamma2}\left|V^{(n)}(t)-\wt{V}^{(n)}(t)\right|&\leq&\vep+4\dl(m_1+1)\sup_{s\in[0,1]}\wt{V}^{(n)}(s)+\omega_{\wt{V}^{(n)}}(6\dl).
\end{eqnarray}

Likewise, we set 
$$\wt{E}_n=E_n\cap\{p_n\geq(1-4\dl(m_1+1))n\}.$$ 
Lemma \ref{Markovsp} implies that, under the probability measure $\ov{\P}^{\,n}_{\mu,\nu,x}(\cdot\mid\wt{E}_n)$ and conditionally on the $\sg$-field $\g_n$ defined by
$$\g_n=\sg\left(\left(\t^{\ov{\al}n^{1/4}},U^{\ov{\al}n^{1/4}}\right),\;M_n,\left(\#\t^{[v_i^n],1},\,1\leq i\leq M_n\right)\right),$$
the spatial tree $(\t^{[v_{i_n}^n]},\ov{U}^{[v_{i_n}^n]})$ is distributed according to $\ov{\P}_{\mu,\nu,Y_n}^{\,p_n}$ (recall that $Y_n=U_{v_{i_n}^n}$). Note that $\wt{E}_n\in\g_n$, and that $Y_n$ and $p_n$ are $\g_n$-measurable. Thus we have
\begin{eqnarray}
\ov{\E}_{\mu,\nu,x}^{\,n}\left(\ind{\wt{E}_n}F\left(\wt{C}^{(n)},\wt{V}^{(n)}\right)\right)&=&\ov{\E}_{\mu,\nu,x}^{\,n}\left(\ind{\wt{E}_n}\ov{\E}_{\mu,\nu,x}^{\,n}\Big(F\left(\wt{C}^{(n)},\wt{V}^{(n)}\right)\Big|\;\g_n\Big)\right)\nonumber\\
&=&\ov{\E}_{\mu,\nu,x}^{\,n}\left(\ind{\wt{E}_n}\ov{\E}_{\mu,\nu,Y_n}^{\,p}\left(F\left(C^{(p)},V^{(p)}\right)\right)_{p=p_n}\right)\label{consMsp}.
\end{eqnarray}
From Lemma \ref{condpos}, we get for every $p$ sufficiently large, 
$$\sup_{\frac{\ov{\al}}{2}p^{1/4}\leq y\leq\frac{3\ov{\al}}{2}p^{1/4}}\left|\ov{\E}_{\mu,\nu,y}^{\,p}\left(F\left(C^{(p)},V^{(p)}\right)\right)-\EE\left(F\left(\ov{\e}^{By/p^{1/4}},\ov{\r}^{By/p^{1/4}}\right)\right)\right|\leq b,$$
which implies using (\ref{contserp}), since $3\ov{\al}B/2\leq2\ov{\al}B=2\al<2\vep$, that for every $p$ sufficiently large,
\be\label{fin1}
\sup_{\frac{\ov{\al}}{2}p^{1/4}\leq y\leq\frac{3\ov{\al}}{2}p^{1/4}}\left|\ov{\E}_{\mu,\nu,y}^{\,p}\left(F\left(C^{(p)},V^{(p)}\right)\right)-\EE\left(F\left(\ov{\e}^{0},\ov{\r}^{0}\right)\right)\right|\leq2b.
\ee
Furthermore Lemma \ref{dperefils} implies that 
$$\ov{\P}_{\mu,\nu,x}^{\,n}\left(\{|n^{-1/4}Y_n-\ov{\al}|>\eta\}\cap E_n\right)\build{\la}_{\to\infty}^{}0.$$ 
So we get for every $n$ sufficiently large, 
\begin{eqnarray}
&&\hspace{-1.3cm}\left|\ov{\E}_{\mu,\nu,x}^{\,n}\left(\ind{\wt{E}_n}\ov{\E}_{\mu,\nu,Y_n}^{\,p}\left(F\left(C^{(p)},V^{(p)}\right)\right)_{p=p_n}\right)-\ov{\P}_{\mu,\nu,x}^{\,n}\left(\wt{E}_n\right)\EE\left(F\left(\ov{\e}^0,\ov{\r}^0\right)\right)\right|\nonumber\\
&\hspace{-0.7cm}\leq&\hspace{-0.5cm}2\ov{\P}_{\mu,\nu,x}^{\,n}\left(\wt{E}_n\cap\left\{|n^{-1/4}Y_n-\ov{\al}|>\ov{\al}/4\right\}\right)\nonumber\\
&&\hspace{-0.2cm}+\;\ov{\E}_{\mu,\nu,x}^{\,n}\left(\ind{\wt{E}_n}\left(\sup_{\frac{\ov{\al}}{2}p^{1/4}\leq y\leq\frac{3\ov{\al}}{2}p^{1/4}}\left|\ov{\E}_{\mu,\nu,y}^{\,p}\left(F\left(C^{(p)},V^{(p)}\right)\right)-\EE\left(F\left(\ov{\e}^{0},\ov{\r}^{0}\right)\right)\right|\right)_{p=p_n}\right)\nonumber\\
&\hspace{-0.7cm}\leq&\hspace{-0.5cm}2b+\ov{\E}_{\mu,\nu,x}^{\,n}\left(\ind{\wt{E}_n}\left(\sup_{\frac{\ov{\al}}{2}p^{1/4}\leq y\leq\frac{3\ov{\al}}{2}p^{1/4}}\left|\ov{\E}_{\mu,\nu,y}^{\,p}\left(F\left(C^{(p)},V^{(p)}\right)\right)-\EE\left(F\left(\ov{\e}^{0},\ov{\r}^{0}\right)\right)\right|\right)_{p=p_n}\right).\label{fin2}
\end{eqnarray}
Thus we use (\ref{consMsp}), (\ref{fin1}), (\ref{fin2}) and the fact that $p_n\geq1-4\dl(m_1+1)n$ on $\wt{E}_n$, to obtain that for every $n$ sufficiently large,
\be\label{fin3}
\left|\ov{\E}_{\mu,\nu,x}^{\,n}\left(\ind{\wt{E}_n}F\left(\wt{C}^{(n)},\wt{V}^{(n)}\right)\right)-\ov{\P}^{\,n}_{\mu,\nu,x}(\wt{E}_n)\EE\left(F\left(\ov{\e}^0,\ov{\r}^0\right)\right)\right|\leq4b.
\ee
From Lemma \ref{p_n}, we have $\ov{\P}_{\mu,\nu,x}^{\,n}(\wt{E}_n)\geq1-2b$. Furthermore, $0\leq F\leq1$ so that (\ref{fin3}) gives
\be\label{fin4}
\left|\ov{\E}_{\mu,\nu,x}^{\,n}\left(F(\wt{C}^{(n)},\wt{V}^{(n)})\right)-\EE(F(\ov{\e}^0,\ov{\r}^0))\right|\leq8b.
\ee
On the other hand, since $\wt{\Gamma}_n\subset\wt{E}_n$ and $F$ is a Lipschitz function whose Lipschitz constant is less than $1$, we have using (\ref{majwtGamma1}) and (\ref{majwtGamma2}), for $n$ sufficiently large, 
\begin{eqnarray}
&&\hspace{-1cm}\ov{\E}_{\mu,\nu,x}^{\,n}\left(\ind{\wt{\Gamma}_n}\left|F\left(\wt{C}^{(n)},\wt{V}^{(n)}\right)-F\left(C^{(n)},V^{(n)}\right)\right|\right)\nonumber\\
&\leq&2\vep+\ov{\E}_{\mu,\nu,x}^{\,n}\left(\left(4\dl(m_1+1)\sup_{s\in[0,1]}\wt{C}^{(n)}(s)\right)\wedge1+\omega_{\wt{C}^{(n)}}(6\dl)\wedge1\right)\nonumber\\
&&\hspace{0.5cm}+\;\ov{\E}_{\mu,\nu,x}^{\,n}\left(\left(4\dl(m_1+1)\sup_{s\in[0,1]}\wt{V}^{(n)}(s)\right)\wedge1+\omega_{\wt{V}^{(n)}}(6\dl)\wedge1\right)\label{fin5}.
\end{eqnarray}
By the same arguments we used to derive (\ref{fin4}), we can bound the right-hand side of (\ref{fin5}), for $n$ sufficiently large, by
\begin{eqnarray*}
&&\hskip-1cm
b+\;2\vep+\;\EE\left(\left(4\dl(m_1+1)\sup_{s\in[0,1]}\ov{\e}^0(s)\right)\wedge1+\omega_{\ov{\e}^0}(6\dl)\wedge1\right)\\
&&\hspace{1cm}+\;\EE\left(\left(4\dl(m_1+1)\sup_{s\in[0,1]}\ov{\r}^0(s)\right)\wedge1+\omega_{\ov{\r}^0}(6\dl)\wedge1\right).
\end{eqnarray*}
From (\ref{veppetit}) together with (\ref{dlm_1}), the latter quantity is bounded above by $7b$. Since 
$$\ov{\P}_{\mu,\nu,x}^{\,n}\left(\wt{\Gamma}_n\right)\geq1-2b,$$ we get for all $n$ sufficiently large,
$$\ov{\E}_{\mu,\nu,x}^{\,n}\left(\left|F\left(\wt{C}^{(n)},\wt{V}^{(n)}\right)-F\left(C^{(n)},V^{(n)}\right)\right|\right)\leq9b,$$
which implies together with (\ref{fin4}) that for all $n$ sufficiently large,
$$\left|\ov{\E}_{\mu,\nu,x}^{\,n}\left(F\left(C^{(n)},V^{(n)}\right)\right)-\EE(F\left(\ov{\e}^0,\ov{\r}^0)\right)\right|\leq17b.$$ 
This completes the proof of Theorem \ref{thconvcond}

\subsection{Proof of Theorem \ref{thcartes}}

In this section we derive Theorem \ref{thcartes} from Theorem \ref{thconvcond}. We first need to prove a lemma. Recall that if $\t\in\T$, we set $\zeta=\#\t-1$ and we denote by $v(0)=\emptyset\prec v(1)\prec\ldots\prec v(\zeta)$ the list of vertices of $\t$ in lexicographical order. For $n\in\{0,1,\ldots,\zeta\}$, we set as in \cite{MaMi},
$$J_\t(n)=\#\left(\t^0\cap\{v(0),v(1),\ldots,v(n)\}\right).$$
We extend $J_\t$ to the real interval $[0,\zeta]$ by setting $J_\t(t)=J_\t(\lfloor t\rfloor)$ for every $t\in[0,\zeta]$, and we set for every $t\in[0,1]$ $$\ov{J}_\t(t)=\frac{J_\t(\zeta t)}{\#\t^0}.$$ 

We also define for $k\in\{0,1,\ldots,2\zeta\}$,
$$K_\t(k)=1+\#\left\{l\in\{1,\ldots,k\}:C(l)=\max_{[l-1,l]}C\;{\rm and}\;C(l)\;{\rm is}\;{\rm even}\right\}.$$
Note that $K_\t(k)$ is the number of vertices of type $0$ in the search-depth sequence up to time $k$. As previously, we extend $K_\t$ to the real interval $[0,2\zeta]$ by setting $K_\t(t)=K_\t(\lfloor t\rfloor)$ for every $t\in[0,2\zeta]$, and we set for every $t\in[0,1]$ 
$$\ov{K}_\t(t)=\frac{K_\t(2\zeta t)}{\#\t^0}.$$ 

\begin{lemma}\label{lemISE}
The law under $\ov{\P}_{\mu,\nu,1}^{\,n}$ of $\left(\ov{J}_\t(t),0\leq t\leq1\right)$ converges as $n\to\infty$ to the Dirac mass at the identity mapping of $[0,1]$. In other words, for every $\eta>0$,
\be\label{eqISEMaMi}
\ov{\P}_{\mu,\nu,1}^{\,n}\left(\sup_{t\in[0,1]}\left|\ov{J}_\t(t)-t\right|>\eta\right)\build{\la}_{n\to\infty}^{}0.
\ee 
Consequently, the law under $\ov{\P}_{\mu,\nu,1}^{\,n}$ of $\left(\ov{K}_\t(t),0\leq t\leq1\right)$ converges as $n\to\infty$ to the Dirac mass at the identity mapping of $[0,1]$. In other words, for every $\eta>0$,
\be\label{eqISEK}
\ov{\P}_{\mu,\nu,1}^{\,n}\left(\sup_{t\in[0,1]}\left|\ov{K}_\t(t)-t\right|>\eta\right)\build{\la}_{n\to\infty}^{}0.
\ee 
\end{lemma}

\proof For $\t\in\T$, we let $v^0(0)=\emptyset\prec v^0(1)\prec\ldots\prec v^0(\#\t^0-1)$ be the list of vertices of $\t$ of type $0$ in lexicographical order. We define as in \cite{MaMi}
$$G_\t(k)=\#\left\{u\in\t:u\prec v^0(k)\right\},\;\;0\leq k\leq\#\t^0-1,$$
and we set $G_\t(\#\t^0)=\zeta$. Note that $v^0(k)$ does not belong to the set $\{u\in\t:u\prec v^0(k)\}$. Recall that $m_0$ denotes the mean of the offspring distribution $\mu_0$. From the second assertion of Lemma 18 in \cite{MaMi} there exists a constant $\vep>0$ such that for $n$ sufficiently large,
$$P_\mu^n\left(\sup_{0\leq k\leq\#\t^0}\left|G_\t(k)-(1+m_0)k\right|\geq n^{3/4}\right)\leq e^{-n^\vep}.$$
Then Lemma \ref{asympt} and Proposition \ref{estpositif} imply that there exists a constant $\vep'>0$ such that for $n$ sufficiently large,
\be\label{eqG_t}
\ov{\P}_{\mu,\nu,1}^{\,n}\left(\sup_{0\leq k\leq\#\t^0}\left|G_\t(k)-(1+m_0)k\right|\geq n^{3/4}\right)\leq e^{-n^{\vep'}}.
\ee
From our definitions, we have for every $0\leq k\leq\#\t^0-1$ and $0\leq n\leq\zeta$,
$$\left\{G_\t(k)>n\right\}=\left\{J_\t(n)\leq k\right\}.$$
It then follows from (\ref{eqG_t}) that, for every $\eta>0$
$$\ov{\P}^{\,n}_{\mu,\nu,1}\left(n^{-1}\sup_{0\leq k\leq\#\t^0}|J_\t(((1+m_0)k)\wedge\zeta)-k|>\eta\right)\build{\la}_{n\to\infty}^{}0.$$
Also from the bound (\ref{beta_1}) of Lemma \ref{feuilles} we get for every $\eta>0$,
$$\ov{\P}^{\,n}_{\mu,\nu,1}\left(\left|n^{-1}\#\t^0-\frac{1}{m_0}\right|>\eta\right)\build{\la}_{n\to\infty}^{}0.$$
The first assertion of Lemma \ref{lemISE} follows from the last two convergences.

Let us set $j_n=2n-|v(n)|$ for $n\in\{0,\ldots,\zeta\}$. It is well known and easy to check by induction that $j_n$ is the first time at which $v(n)$ appears in the search-depth sequence. It is also convenient to set $j_{\zeta+1}=2\zeta$. Then we have $K_\t(k)=J_\t(n)$ for every $k\in\{j_n,\ldots,j_{n+1}-1\}$ and every $n\in\{0,\ldots,\zeta\}$. Let us define a random function $\vph:[0,2\zeta]\la\Z_+$ by setting $\vph(t)=n$ if $t\in[j_n,j_{n+1})$ and $0\leq n\leq\zeta$, and $\vph(2\zeta)=\zeta$. From our definitions, we have for every $t\in[0,2\zeta]$,
\be\label{eqISE2}
K_\t(t)=J_\t(\vph(t)).
\ee 
Furthermore, we easily check from the equality $j_n=2n-|v(n)|$ that
\be\label{eqISE1}
\sup_{t\in[0,2\zeta]}\left|\vph(t)-\frac{t}{2}\right|\leq\max_{[0,2\zeta]}C.
\ee
We set $\vph_\zeta(t)=\zeta^{-1}\vph(2\zeta t)$ for $t\in[0,1]$. So (\ref{eqISE1}) gives
$$\sup_{t\in[0,1]}|\vph_\zeta(t)-t|\leq\frac{1}{\zeta}\,\max_{[0,2\zeta]}C,$$
which implies that for every $\eta>0$,
\be\label{eqISE5}
\ov{\P}_{\mu,\nu,1}^{\,n}\left(\sup_{t\in[0,1]}|\vph_\zeta(t)-t|>\eta\right)\build{\la}_{n\to\infty}^{}0.
\ee 
On the other hand we get from (\ref{eqISE2}) that $\ov{K}_\t(t)=\ov{J}_\t(\vph_\zeta(t))$ for every $t\in[0,1]$ and thus
$$\sup_{t\in[0,1]}\left|\ov{K}_\t(t)-\ov{J}_\t(t)\right|\leq2\sup_{t\in[0,1]}\left|\ov{J}_\t(t)-t\right|+\sup_{t\in[0,1]}|\vph_\zeta(t)-t|.$$
The desired result then follows from (\ref{eqISEMaMi}) and (\ref{eqISE5}).\cq
 
We are now able to complete the proof of Theorem \ref{thcartes}. The proof of (i) is similar to the proof of the first part of Theorem 8.2 in \cite{LG}, and is therefore omitted.  

Let us turn to (ii). By Corollary \ref{loiimage} and properties of the Bouttier-di Francesco-Guitter bijection, the law of $\lm_M^{(n)}$ under $\B_\q^r(\cdot\mid\#\f_M=n)$ is the law under $\ov{\P}_{\mu,\nu,1}^{\,n}$ of the probability measure $\I_n$ defined by
$$\langle\I_n,g\rangle=\frac{1}{\#\t^0+1}\left(g(0)+\sum_{v\in\t^0}g\left(n^{-1/4}U_v\right)\right).$$
It is more convenient for our purposes to replace $\I_n$ by a new probability measure $\I'_n$ defined by
$$\langle\I'_n,g\rangle=\frac{1}{\#\t^0}\sum_{v\in\t^0}g\left(n^{-1/4}U_v\right).$$
Let $g$ be a bounded continuous function. Clearly, we have for every $\eta>0$, 
\be\label{eqconvISE1}
\ov{\P}_{\mu,\nu,1}^{\,n}\left(\left|\langle\I_n,g\rangle-\langle\I'_n,g\rangle\right|>\eta\right)\build{\la}_{n\to\infty}^{}0.
\ee 
Furthermore, we have from our definitions
\be\label{defISE}
\langle\I'_n,g\rangle=\frac{1}{\#\t^0}\,g\left(n^{-1/4}\right)+\int_0^1g\left(n^{-1/4}V(2\zeta t)\right)\d\ov{K}_\t(t),
\ee
where the first term in the right-hand side corresponds to $v=\emptyset$ in the definition of $\I'_n$. Then from Theorem \ref{thconvcond}, (\ref{eqISEK}) and the Skorokhod representation theorem, we can construct on a suitable probability space, a sequence $(\t_n,U_n)_{n\geq1}$ and a conditioned Brownian snake $(\ov{\eb}^0,\ov{\rb}^0)$, such that each pair $(\t_n,U_n)$ is distributed according to $\ov{\P}_{\mu,\nu,1}^{\,n}$, and such that if we write $(C_n,V_n)$ for the  contour functions of $(\t_n,U_n)$, $\zeta_n=\#\t_n-1$ and $\ov{K}_n=\ov{K}_{\t_n}$, we have,
$$\left(\frac{V_n(2\zeta_nt)}{n^{1/4}},\ov{K}_n(t)\right)\build{\la}_{n\to\infty}^{}\left(\left(\frac{4\rho_\q}{9(Z_\q-1)}\right)^{1/4}\,\ov{\rb}^0(t),\;t\right),$$
uniformly in $t\in[0,1]$, a.s. Now $g$ is Lipschitz, which implies that a.s.,
\be\label{eqconvISE2}
\left|\int_0^1g\left(n^{-1/4}V_n(2\zeta_nt)\right)\d\ov{K}_n(t)-\int_0^1g\left(\left(\frac{4\rho_\q}{9(Z_\q-1)}\right)^{1/4}\,\ov{\rb}^0(t)\right)\d\ov{K}_n(t)\right|\build{\la}_{n\to\infty}^{}0.
\ee
Furthermore, the sequence of measures $\d\ov{K}_n$ converges weakly to the uniform measure $\d t$ on $[0,1]$ a.s., so that a.s.
\be\label{eqconvISE3}
\int_0^1g\left(\left(\frac{4\rho_\q}{9(Z_\q-1)}\right)^{1/4}\,\ov{\rb}^0(t)\right)\d\ov{K}_n(t)\build{\la}_{n\to\infty}^{}\int_0^1g\left(\left(\frac{4\rho_\q}{9(Z_\q-1)}\right)^{1/4}\,\ov{\rb}^0(t)\right)\d t.
\ee
Then (\ref{eqconvISE2}) and (\ref{eqconvISE3}) imply that a.s.,
$$\int_0^1g\left(n^{-1/4}V_n(2\zeta_nt)\right)\d\ov{K}_n(t)\build{\la}_{n\to\infty}^{}\int_0^1g\left(\left(\frac{4\rho_\q}{9(Z_\q-1)}\right)^{1/4}\,\ov{\rb}^0(t)\right)\d t,$$
which together with (\ref{eqconvISE1}) and (\ref{defISE}) yields the desired result. 

Finally, the proof of (iii) from (ii) is similar to the proof of the third part of Theorem 8.2 in \cite{LG}. This completes the proof of Theorem \ref{thcartes}.

\section{Separating vertices in a $2\ka$-angulation}\label{secdisc}

In this section, we use the estimates of Proposition \ref{estpositif} to derive a result concerning separating vertices in rooted $2\ka$-angulations. Recall that in a $2\ka$-angulation, all faces have a degree equal to $2\ka$. 

Let $M$ be a planar map and let $\sg_0\in\v_M$. Let $\sg$ be a vertex of $M$ different from $\sg_0$. We denote by $\s_M^{\sg_0,\sg}$ the set of all vertices $a$ of $M$ such that any path from $\sg$ to $a$ goes through $\sg_0$. The vertex $\sg_0$ is called a {\em separating vertices} of $M$ if there exists a vertex $\sg$ of $M$ different from $\sg_0$ such that $\s_M^{\sg_0,\sg}\neq\{\sg_0\}$. We denote by $\DD_M$ the set of all separating vertices of $M$.

Recall that $\ov{\U}_{\ka}^{\,n}$ stands for the uniform probability measure on the set of all rooted $2\ka$-angulations with $n$ faces. Our goal is to prove the following theorem.

\begin{theorem}\label{disccartesm_r}
For every $\vep>0$,
$$\lim_{n\to\infty}\ov{\U}_\ka^{\,n}\left(\exists\,\sg_0\in\DD_M:\exists\,\sg\in\v_M\setminus\{\sg_0\},\,n^{1/2-\vep}\leq\#\s_M^{\sg_0,\sg}\leq2n^{1/2-\vep}\right)=1.$$
\end{theorem}

Theorem \ref{disccartesm_r} is a consequence of the following theorem. Recall that $\U_{\ka}^{n}$ denotes the uniform probability measure on the set of all rooted pointed $2\ka$-angulations with $n$ faces. If $M$ is a rooted pointed bipartite planar map, we denote by $\tau$ its distinguished point.

\begin{theorem}\label{disccartes}
For every $\vep>0$,
$$\lim_{n\to\infty}\U_\ka^{n}\left(\exists\,\sg_0\in\DD_M:\sg_0\neq\tau,\,n^{1/2-\vep}\leq\#\s_M^{\sg_0,\tau}\leq2n^{1/2-\vep}\right)=1.$$
\end{theorem}

Theorem \ref{disccartesm_r} can be deduced from Theorem \ref{disccartes} but not as directly as one could think. Indeed the canonical surjection from the set of rooted pointed $2\ka$-angulations with $n$ faces onto the set of rooted $2\ka$-angulations with $n$ faces does not map the uniform measure $\U^n_\ka$ to the uniform measure $\ov{\U}^{\,n}_\ka$. Nevertheless a simple argument allows us to circumvent this difficulty. 

Let $\wt{\m}_{r,p}$ be the set of all triples $(M,\vec{e},\tau)$ where $(M,\vec{e}\,)\in\m_r$ and $\tau$ is a distinguished vertex of the map $M$. We denote by $\ss$ the canonical surjection from the set $\wt{\m}_{r,p}$ onto the set $\m_{r,p}$ which is obtained by ``forgetting'' the orientation of $\vec{e}$. We observe that for every $(M,e,\tau)\in\m_{r,p}$ 
$$\#\left(\ss^{-1}((M,e,\tau))\right)=2.$$
Denote by $\wt{\U}^{\,n}_\ka$ the uniform measure on the set of all triples $(M,\vec{e},\tau)\in\wt{\m}_{r,p}$ such that $M$ is a $2\ka$-angulation with $n$ faces. Then the image measure of the measure $\wt{\U}^{\,n}_{\ka}$ under the mapping $\ss$ is the measure $\U_\ka^n$. Thus we obtain from Theorem \ref{disccartes} that
\be\label{eqthdisccartesm_r}
\lim_{n\to\infty}\wt{\U}^{\,n}_\ka\left(\exists\,\sg_0\in\DD_M:\exists\,\sg\in\v_M\setminus\{\sg_0\},\,n^{1/2-\vep}\leq\#\s_M^{\sg_0,\sg}\leq2n^{1/2-\vep}\right)=1.
\ee
On the other hand let $\pp$ be the canonical projection from the set $\wt{\m}_{r,p}$ onto the set $\m_{r}$. If $M$ is a $2\ka$-angulation with $n$ faces, we have thanks to Euler formula  
$$\#\v_M=(\ka-1)n+2.$$
Thus the image measure of the measure $\wt{\U}^{\,n}_{\ka}$ under the mapping $\pp$ is the measure $\ov{\U}_\ka^{\,n}$. This remark together with (\ref{eqthdisccartesm_r}) implies Theorem \ref{disccartesm_r}.



The remainder of this section is devoted to the proof of Theorem \ref{disccartes}. We first need to state a lemma. Recall the definition of the spatial tree $(\t^{[v]},U^{[v]})$ for $(\t,U)\in\Omega$ and $v\in\t$.

\begin{lemma}\label{discarbres}
Let $(\t,U)\in\T_1^{\rm mob}$ and let $M=\Psi_{r,p}((\t,U))$. Suppose that we can find $v\in\t^0$ such that $\t^{[v],0}\neq\{\emptyset\}$ and $\un{U}^{[v]}>0$. Then there exists $\sg_0\in\DD_M$ such that $\sg_0\neq\tau$ and
$$\#\s_M^{\sg_0,\tau}=\#\t^{[v],0}.$$
\end{lemma}

\proof Let $(\t,U)\in\T_1^{\rm mob}$. Write $w_0,w_1,\ldots,w_\zeta$ for the search-depth sequence of $\t^{0}$ (see Section \ref{secbij}). Recall from Section \ref{secbij} the definition of $(U^+_v,v\in\t)$ and the construction of the planar map $\Psi_{r,p}((\t,U))$. For every $i\in\{0,1,\ldots,\zeta\}$, we set $s_i=\dd$ if $U^+_{w_i}=1$, whereas if $U^+_{w_i}\geq2$, we denote by $s_i$ the first vertex in the sequence $w_{i+1},\ldots,w_{\zeta-1},w_0,w_1,\ldots,w_{i-1}$ whose label is $U^+_{w_i}-1$. 

Suppose that there exists $v\in\t^0$ such that $\t^{[v],0}\neq\{\emptyset\}$ and $\un{U}^{[v]}>0$. We set
\ba
k&=&\min\{i\in\{0,1,\ldots,\zeta\}:w_i=v\},\\
l&=&\max\{i\in\{0,1,\ldots,\zeta\}:w_i=v\}.
\ea
The vertices $w_k,w_{k+1},\ldots,w_l$ are exactly the descendants of $v$ in $\t^0$. The condition $\un{U}^{[v]}>0$ ensures that for every $i\in\{k+1,\ldots,l-1\}$, we have
$$U^+_{w_i}>U^+_{w_l}=U^+_{w_k}.$$
This implies that $s_i$ is a descendant of $v$ for every $i\in\{k+1,\ldots,l-1\}$. Furthermore $s_k=s_l$ and $s_i$ is not a strict descendant of $v$ if $i\in\{0,1,\ldots,\zeta\}\setminus\{k,k+1\ldots,l\}$. From the construction of edges in the map $\Psi_{r,p}((\t,U))$ we see that any path from $\dd$ to a vertex that is a descendant of $v$ must go through $v$. It follows that $v$ is a separating vertex of the map $M=\Psi_{r,p}((\t,U))$ and that the set $\t^{[v],0}$ is in one-to-one correspondence with the set $S^{v,\tau}_M$.  \cq

Thanks to Corollary \ref{loiimageka} and Lemma \ref{discarbres}, it suffices to prove the following proposition in order to get Theorem \ref{disccartes}. Recall the definition of $\mu^\ka=(\mu_0^\ka,\mu_1^\ka)$.

\begin{proposition}\label{propdiscarbres}
For every $\vep>0$,
$$\lim_{n\to\infty}\P_{\mu^\ka,\nu,1}^{n}\left(\exists\,v_0\in\t^0:n^{1/2-\vep}\leq\#\t^{[v_0],0}\leq2n^{1/2-\vep},\;\un{U}^{[v_0]}>0\right)=1.$$
\end{proposition}

\proof For $n\geq1$ and $\vep>0$, we denote by $\Lambda_{n,\vep}$ the event
$$\Lambda_{n,\vep}=\left\{\exists\,v_0\in\t^0:n^{1/2-\vep}\leq\#\t^{[v_0],1}\leq2n^{1/2-\vep},\;\un{U}^{[v_0]}>0\right\}.$$
Let $\vep>0$ and $\al>0$. We will prove that for all $n$ sufficiently large,
\be\label{eqdisc0}
\P_{\mu^\ka,\nu,1}^n(\Lambda_{n,\vep})\geq1-3\al.
\ee
We first state a lemma. For $\t\in\T$ and $k\geq0$, we set 
$$Z^0(k)=\#\{v\in\t:|v|=2k\}.$$

\begin{lemma}\label{lemdisc1}
There exist constants $\beta>0$, $\gamma>0$ and $M>0$ such that for all $n$ sufficiently large,
$$P^{n}_{\mu^\ka}\left(\inf_{\gamma\sqrt{n}\leq k\leq2\gamma\sqrt{n}}Z^0(k)>\beta\sqrt{n},\;\;\sup_{k\geq0}Z^0(k)<M\sqrt{n}\right)\geq1-\al.$$
\end{lemma}

We postpone the proof of Lemma \ref{lemdisc1} and complete that of Proposition \ref{propdiscarbres}. To this end, we introduce some notation. For $n\geq1$, we define an integer $K_n$ by the condition
$$\lceil\gamma\sqrt{n}\,\rceil+K_n\lceil n^{1/4}\rceil\leq\lfloor2\gamma\sqrt{n}\rfloor<\lceil\gamma\sqrt{n}\,\rceil+(K_n+1)\lceil n^{1/4}\rceil,$$
and we set for $j\in\{0,\ldots,K_n\}$,
$$k^{(n)}_j=\lceil\gamma\sqrt{n}\,\rceil+j\lceil n^{1/4}\rceil.$$
If $\t\in\T$, we write $\h(\t)$ for the height of $\t$, that is the maximal generation of an individual in $\t$. For $k\geq0$, and $N,P\geq1$ we set
$$Z^0(k,N,P)=\#\left\{v\in\t:|v|=2k,N\leq\#\t^{[v],0}\leq2N,\h(\t^{[v]})\leq2P\right\}.$$
We denote by $\Gamma_n$, $C_{n,\vep}$ and $E_{n,\vep}$ the events
\ba
\Gamma_n&=&\bigcap_{0\leq j\leq K_n}\left\{\beta\sqrt{n}<Z^0\left(k_j^{(n)}\right)<M\sqrt{n}\right\},\\
C_{n,\vep}&=&\bigcap_{0\leq j\leq K_n}\left\{Z^0\left(k_j^{(n)},n^{1/2-\vep},n^{1/4}\right)>n^{1/4}\right\},\\
E_{n,\vep}&=&\bigcap_{0\leq j\leq K_n}\left\{n^{1/4}<Z^0\left(k_j^{(n)},n^{1/2-\vep},n^{1/4}\right)<M\sqrt{n}\right\}.
\ea
The first step is to prove that for all $n$ sufficiently large,
\be\label{eqdisc1}
P_{\mu^\ka}^{\,n}\left(E_{n,\vep}\right)\geq1-2\al.
\ee
Since $\Gamma_n\cap C_{n,\vep}\subset E_{n,\vep}$, it suffices to prove that for all $n$ sufficiently large,
\be\label{eqdisc2}
P_{\mu^\ka}^{\,n}\left(\Gamma_n\cap C_{n,\vep}\right)\geq1-2\al.
\ee
We first observe that
\be\label{eqdisc3}
P_{\mu^\ka}(\Gamma_n\cap C_{n,\vep}^c)\leq\sum_{j=0}^{K_n}P_{\mu^\ka}\left(\beta\sqrt{n}<Z^0\left(k_j^{(n)}\right)<M\sqrt{n},\;Z^0\left(k_j^{(n)},n^{1/2-\vep},n^{1/4}\right)\leq n^{1/4}\right).
\ee
Let $j\in\{0,\ldots,K_n\}$. We have
\begin{eqnarray}
&&\hspace{-1.2cm}P_{\mu^\ka}\left(\beta\sqrt{n}<Z^0\left(k_j^{(n)}\right)<M\sqrt{n},\;Z^0\left(k_j^{(n)},n^{1/2-\vep},n^{1/4}\right)\leq n^{1/4}\right)\nonumber\\
&\hspace{-0.5cm}=&\hspace{-0.5cm}\sum_{q=\lceil\beta\sqrt{n}\,\rceil}^{\lfloor M\sqrt{n}\,\rfloor}P_{\mu^\ka}\left(Z\left(k_j^{(n)},n^{1/2-\vep},n^{1/4}\right)\leq n^{1/4}\;\;\Big|\;\,Z^0\left(k_j^{(n)}\right)=q\right)P_{\mu^\ka}\left(Z^0\left(k_j^{(n)}\right)=q\right).\label{eqdisc4}
\end{eqnarray}
Now, under the probability measure $P_{\mu^\ka}(\cdot\mid Z^0(k_j^{(n)})=q)$, the $q$ subtrees of $\t$ above level $2k_j^{(n)}$ are independent and distributed according to $P_{\mu^\ka}$. Consider on a probability space $(\Omega',\PP')$ a sequence of Bernoulli variables $(B_i,i\geq1)$ with parameter $p_n$ defined by
$$p_n=P_{\mu^\ka}\left(n^{1/2-\vep}\leq\#\t^0\leq2n^{1/2-\vep},\h(\t)\leq2n^{1/4}\right).$$
Then (\ref{eqdisc4}) gives
\begin{eqnarray}
&&\hspace{-1.5cm}P_{\mu^\ka}\left(\beta\sqrt{n}<Z^0\left(k_j^{(n)}\right)<M\sqrt{n},\;Z\left(k_j^{(n)},n^{1/2-\vep},n^{1/4}\right)\leq n^{1/4}\right)\nonumber\\
&\leq&\sum_{q=\lceil\beta\sqrt{n}\,\rceil}^{\lfloor M\sqrt{n}\,\rfloor}\PP'\left(\sum_{i=1}^qB_i\leq n^{1/4}\right)\nonumber\\
&\leq&\sum_{q=\lceil\beta\sqrt{n}\,\rceil}^{\lfloor M\sqrt{n}\,\rfloor}e\exp\left(-\frac{qp_nn^{-1/4}}{2}\right)\label{eqdisc5},
\end{eqnarray}
where the last bound follows from a simple exponential inequality. However, from Lemma 14 in \cite{MaMi}, there exists $\eta>0$ such that for all $n$ sufficiently large,
\begin{eqnarray}
&\hspace{0.8cm}p_n&=\;P_{\mu^\ka}\left(n^{1/2-\vep}\leq\#\t^0\leq2n^{1/2-\vep}\right)-P_\mu\left(n^{1/2-\vep}\leq\#\t^0\leq2n^{1/2-\vep},\h(\t)>2n^{1/4}\right)\nonumber\\
&&\geq\;P_{\mu^\ka}\left(n^{1/2-\vep}\leq\#\t^0\leq2n^{1/2-\vep}\right)-e^{-n^\eta}\label{MaMi14}.
\end{eqnarray}
Under the probability measure $P_{\mu^\ka}$, we have $\#\t^0=(\ka-1)\#\t^1+1$ a.s., so we get from Lemma \ref{asympt} that there exists a constant $c_\ka$ such that,
\be\label{eqasymptdisc}
n^{1/4-\vep/2}\,P_{\mu^\ka}\left(n^{1/2-\vep}\leq\#\t^0\leq2n^{1/2-\vep}\right)\build{\la}_{n\to\infty}^{}c_{\ka}.
\ee
It then follows from (\ref{MaMi14}) that
$$n^{1/4-\vep/2}\,p_n\build{\la}_{n\to\infty}^{}c_\ka.$$
From (\ref{eqdisc5}), we obtain that there exists a constant $c'>0$ such that 
$$P_{\mu^\ka}\left(\beta\sqrt{n}<Z^0_{k_j^{(n)}}<M\sqrt{n},\;Z\left(k_j^{(n)},n^{1/2-\vep},n^{1/4}\right)\leq n^{1/4}\right)\leq M\sqrt{n}e^{-c'n^\vep},$$
which together with (\ref{eqdisc3}) implies that 
$$P_{\mu^\ka}\left(\Gamma_n\cap C^c_{n,\vep}\right)\leq MK_n\sqrt{n}e^{-n^\vep}.$$ 
Since $K_n\sim\gamma n^{1/4}$ as $n\to\infty$, we get from Lemma \ref{asympt} that, for all $n$ sufficiently large,
\be\label{eqdisc6}
P_{\mu^\ka}^n\left(\Gamma_n\cap C_{n,\vep}^c\right)\leq\frac{P_{\mu^\ka}(\Gamma_n\cap C_{n,\vep}^c)}{P_{\mu^\ka}(\#\t^1=n)}\leq\al.
\ee
On the other hand, Lemma \ref{lemdisc1} implies that 
$$P_{\mu^\ka}^{\,n}(\Gamma_n)\geq1-\al.$$ 
The bound (\ref{eqdisc2}) now follows.

Set 
$$I_n=\left\{\lceil n^{1/4}\,\rceil,\lceil n^{1/4}\,\rceil+1,\ldots,\lfloor M\sqrt{n}\rfloor\right\}^{K_n+1}.$$ 
For $(p_0,\ldots,p_{K_n})\in I_n$, we set 
$$E_{n,\vep}(p_0,\ldots,p_{K_n})=\bigcap_{0\leq j\leq K_n}\left\{Z^0\left(k_j^{(n)},n^{1/2-\vep},n^{1/4}\right)=p_j\right\}.$$
On the event $E_{n,\vep}(p_0,\ldots,p_{K_n})$ for every $j\in\{1,\ldots,K_n\}$, let $v^j_1\prec\ldots\prec v^j_{p_j}$ be the list in lexicographical order of those vertices $v\in\t$ at generation $2k_j^{(n)}$ such that $n^{1/2-\vep}\leq\#\t^{[v],0}\leq2n^{1/2-\vep}$ and $\h(\t^{[v]})\leq2n^{1/4}$. Note that for every $j,j'\in\{1,\ldots,K_n\}$ such that $j<j'$, we have for every $i\in\{1,\ldots,p_j\}$,
$$\max\left\{|v_i^jv|:v\in\t^{[v_i^j]}\right\}\leq2k_j^{(n)}+2n^{1/4}<2k_{j'}^{(n)}.$$
Then, it is not difficult to check that under $\P_{\mu^\ka,\nu,1}(\cdot\mid E_{n,\vep}(p_0,\ldots,p_{K_n}))$, the spatial trees $\{(\t^{[v^j_i]},U^{[v^j_i]}),\,i=1,\ldots,p_j,\,j=1,\ldots,K_n\}$ are independent and distributed according to the probability measure $\P_{\mu^\ka,\nu,0}(\cdot\mid n^{1/2-\vep}\leq\#\t^0\leq n^{1/2-\vep},\h(\t)\leq2n^{1/4})$. Set
$$\pi_n=\P_{\mu,\nu,0}\left(\un{U}>0\mid n^{1/2-\vep}\leq\#\t^0\leq 2n^{1/2-\vep},\h(\t)\leq2n^{1/4}\right).$$
Since the events $E_{n,\vep}(p_0,\ldots,p_{K_n}),\;(p_0,\ldots,p_{K_n})\in I_n$ are disjoints we have,
\begin{eqnarray}
\P_{\mu^\ka,\nu,1}(E_{n,\vep}\cap\Lambda_{n,\vep}^c)&=&\sum_{(p_0,\ldots,p_{K_n})\in I_n}P_{\mu^\ka}(E_{n,\vep}(p_0,\ldots,p_{K_n}))\nonumber\\
&&\times\;\P_{\mu^\ka,\nu,1}\left(\left\{\un{U}^{[v^j_i]}\leq0:0\leq j\leq K_n;1\leq i\leq p_j\right\}\mid E_{n,\vep}(p_0,\ldots,p_{K_n})\right)\nonumber\\
&=&\sum_{(p_0,\ldots,p_{K_n})\in I_n}P_{\mu^\ka}(E_{n,\vep}(p_0,\ldots,p_{K_n}))(1-\pi_n)^{p_0+\ldots+p_{K_n}}\nonumber\\
&\leq&P_{\mu^\ka}(E_{n,\vep})(1-\pi_n)^{\lceil n^{1/4}\rceil(K_n+1)}\nonumber\\
&\leq&(1-\pi_n)^{\lceil n^{1/4}\rceil(K_n+1)}\label{eqdisc7}.
\end{eqnarray}
Now, from Lemma 14 in \cite{MaMi} and $(\ref{eqasymptdisc})$ there exists $c''>0$ such that for all $n$ sufficiently large,
$$\pi_n\geq\P_{\mu,\nu,0}\left(\un{U}>0\mid n^{1/2-\vep}\leq\#\t^0\leq 2n^{1/2-\vep}\right)-c''n^{1/2}e^{-n^\eta},$$
where $\eta$ was introduced before (\ref{MaMi14}). Since under the probability measure $P_{\mu^\ka}$, we have $\#\t^0=1+(\ka-1)\#\t^1$ a.s., we get from Proposition \ref{estpositif} that there exists a constant $\un{c}>0$ such that for all $n$ sufficiently large, 
$$\pi_n\geq\frac{\un{c}}{n^{1/2-\vep}}.$$
So, there exists a constant $\un{c}'>0$ such that (\ref{eqdisc7}) becomes for all $n$ sufficiently large, 
\be\label{eqdisc8}
\P_{\mu^\ka,\nu,1}\left(E_{n,\vep}\cap\Lambda_{n,\vep}^c\right)\leq e^{-\un{c}'n^{\vep}}.
\ee
Finally (\ref{eqdisc8}) together with Lemma \ref{asympt} implies that if $n$ is sufficiently large,
$$\P_{\mu^\ka,\nu,1}^{\,n}\left(E_{n,\vep}\cap\Lambda_{n,\vep}^c\right)\leq\al.$$
Using also (\ref{eqdisc1}), we obtain our claim (\ref{eqdisc0}). \cq

\noindent\textsc{Proof of Lemma \ref{lemdisc1}~:} Under $P_{\mu^\ka}$, $Z^0=(Z^0(k),k\geq0)$ is a critical Galton-Watson process with offspring law $\mu^\ka_{0,\ka}$ supported on $(\ka-1)\Z_+$ and defined by 
$$\mu^\ka_{0,\ka}((\ka-1)k)=\mu^\ka_0(k),\;\;k\geq0.$$ 
Note that the total progeny of $Z^0$ is $\#\t^0$ and recall that under the probability measure $P_{\mu^\ka}$, we have $\#\t^0=1+(\ka-1)\#\t^1$ a.s.

Define a function $L=(L(t),t\geq0)$ by interpolating linearly $Z^0$ between successive integers. For $n\geq1$ and $t\geq0$, we set
$$l_n(t)=\frac{1}{\sqrt{n}}L(t\sqrt{n}).$$
Denote by $\lg(t)$ the total local time of $\e$ at level $t$ that is,
$$\lg(t)=\lim_{\vep\to0}\,\frac{1}{\vep}\int_0^1\ind{\{t\leq\e(s)\leq t+\vep\}}\,\d s,$$
where the convergence holds a.s. From Theorem 1.1 in \cite{DG}, the law of $(l_n(t),t\geq0)$ under $P_{\mu^\ka}^{\,n}$ converges to the law of $(\lg(t),t\geq0)$. So we have
$$\liminf_{n\to\infty}P_{\mu^\ka}^{\,n}\left(\inf_{\gamma\leq t\leq2\gamma}l_n(t)>\beta,\;\;\sup_{t\geq0}\,l_n(t)<M\right)\geq\PP\left(\inf_{\gamma\leq t\leq2\gamma}\lg(t)>\beta,\;\;\sup_{t\geq0}\,\lg(t)<M\right).$$
However we can find $\beta>0$, $\gamma>0$ and $M>0$ such that 
$$\PP\left(\inf_{\gamma\leq t\leq2\gamma}\lg(t)>\beta,\;\;\sup_{t\geq0}\,\lg(t)<M\right)\geq1-\frac{\al}{2},$$
and the desired result follows.\cq

\end{document}